\documentclass[leqno]{article}
\usepackage{graphicx}
\usepackage{booktabs}
\usepackage{amsfonts}
\usepackage{amssymb}
\usepackage{amsthm}
\usepackage{mathrsfs}
\usepackage{amsfonts}
\usepackage{enumerate}
\usepackage{xcolor}
\usepackage{graphicx}
\usepackage{graphics}
\usepackage{amsfonts}
\usepackage{amssymb}
\usepackage{amsthm}
\usepackage{mathrsfs}
\usepackage{amsfonts}
\usepackage{enumerate}
\usepackage{blindtext}
\usepackage{amsmath}
\usepackage{latexsym}
\usepackage{graphics}
\usepackage{graphicx}
\usepackage{wrapfig}
\usepackage{cancel}
\usepackage{ulem}
\usepackage{soul}
\usepackage[polutonikogreek,english]{babel}
\usepackage{scalerel}
\numberwithin{equation}{section}

\newcommand{\ad}{\begin{description}}
\newcommand{\zd}{\end{description}}

\theoremstyle{definition}

\def\widen{\mathrel{\ThisStyle{\stretchrel*{\ooalign{%
  \raise0.6\LMex\hbox{$\SavedStyle\sqcup$}\cr%
  \raise-0.2\LMex\hbox{$\SavedStyle\sqcup$}}}{\sqcup}}}}
\begin{document}

\title{Frege's Theory of Real Numbers: A Consistent Rendering}

\author{FRANCESCA BOCCUNI\footnote{Vita-Salute San Raffaele University, Milan. Orcid: 0000-0001-9814-1431.} ~\& MARCO PANZA\footnote{IHPST (CNRS and Universit\'{e} de Paris 1, Panth\'{e}one Sorbonne), and Chapman University, Orange, CA.}}
\date{}

\maketitle

\begin{abstract}
\noindent Frege's definition of the real numbers, as envisaged in the second volume of \textit{Grundgesetze der Arithmetik}, is fatally flawed by the inconsistency of Frege's ill-fated \textit{Basic Law V}. We restate Frege's definition in a consistent logical framework and investigate whether it can provide a logical foundation of real analysis. Our conclusion will deem it doubtful that such a foundation along the lines of Frege's own indications is possible at all.
\end{abstract}

\section{Overview}
The aim of the present paper is twofold: (\textit{i})~rephrasing Frege's inconsistent definition\footnote{In what follows we will
use `definition' quite broadly. The meaning of its occurrences will be clarified in context.} of real numbers, as envisaged in
Part III of \textit{Grundgesetze der Arithmetik} (Frege 1893-1903), in a consistent setting ruling out value-ranges, and so involving 
no version of the infamous \textit{Basic Law V} (BLV); (\textit{ii})~wondering whether the rephrased definition can be considered logical, and, as such, as a ground for a logicist view about real analysis. 

Concerning (\textit{ii}) a proviso is in order. In the debate on neologicism, a distinction has been made between logicality and analyticity, by suggesting, for instance, that, though not logical, Hume's Principle (HP) is analytic. 
We are far from undermining the relevance of this distinction, but we consider unnecessary to stress it for our present purpose. 
There are two reasons for that. On the one side, we deem
all the arguments we will advance against the logicality of the relevant principles and definitions
also apt to oppose their analyticity---though some of those advanced in favor of
the former are possibly only sufficient to argue for the latter.  
On the other side, we are interested in the epistemic attitude that a faithful Fregean (or even Frege himself) might have
(had) in the face of a definition such as our own. Hence, for the sake of our discussion, we must follow Frege himself in
taking a ``truth'' to be 
analytic if, in its proof, ``one only runs
into logical laws and definitions'' (Frege 1884, \S~3; Frege 1953, p.~4)\footnote{We slightly modify Austin's 
translation.}, and in regarding definitions as mandatorily explicit, which suggests regarding logicality as a
necessary condition for analyticity, rather than the latter as a weaker condition than the former.

Concerning (\textit{i}), it is important to observe that, for Frege, real numbers had 
to be defined as ratios of magnitudes, and
magnitudes had to belong to different domains. Hence,
his definition should have included two successive steps: a
definition of domains of magnitudes, and a
definition of ratios on these domains.
As a matter of fact, he
accomplished only the former step, and merely gave some
informal indications on how to accomplish the latter. 
Both things are done in the
second volume of \textit{Grundgesetze}\footnote{Frege (1893-1903),~\S\S~II.165-II.245  and \S~II.164, respectively; \S\S~II.55-II.159 contain a critical discussion of alternative definitions, while \S\S~II.160-II.163 contain an informal
introduction and a principled justification of the definition of domains of magnitudes.}.
The latter step should have presumably been accomplished in a third volume that, once aware of 
Russell's paradox, Frege never wrote.

Had he accomplished this step, he should have made it conform to a crucial requirement: supposing that 
several domains of magnitude exist, the definition of real
numbers should have identified a ratio on one of
these domains with the same real number as a ratio on each one of the others. Our definition
actually complies with this requirement.

After a short presentation of Frege's
strategy in \S~2,
we will consistently rephrase his definition of
domains of magnitudes in \S~3\footnote{Knowledge of
Frege's original definition is required to appreciate its correspondence
with our rephrasing. Useful accounts of it can be found in
Dummett (1991, ch.~22); Schirn (2013); Simons (1987); Shapiro \& Snyder (2020).}. 
To eliminate value-ranges, we will rephrase first-order formulas involving terms for them as
higher-order formulas proper to a system of higher-order predicate logic as weak as
possible, on which we will take stock in \S~4. 
This seems to us the most faithful way to consistently render Frege's original definition. 
Insofar as our appreciation of the logicality of our definition
depends on
assuming, in a genuinely Fregean vein, the logicality of higher-order logic, we contend that
this appreciation \textit{ipso facto} provides an appreciation of the logicality of Frege's own
definition that remains perfectly independent of any judgement about the
logicality of (any consistent version of) BLV\footnote{Anyone supporting Quine's view on 
the non-logicality of higher-order logic can take our granting it as made for the sake of the argument.\label{Quine}}.

In \S~5, we will investigate how to define real numbers by following 
Frege's indications, on the base of our definition of domains of magnitudes. In
\S~II.164 of his treatise\label{164}, Frege explicitly acknowledged that his
envisaged definition of these numbers requires an existence proof of nonempty such domains. We will explain why this is so.
Here, it is only in order to observe that, in this same \S, he also argues that
this existence depends on the existence of continuously many objects (an infinity of objects larger than ``\textit{Endloss}'', the
cardinality of ``finite cardinal numbers''), and
sketches a plan for this proof, which, taking the existence of natural
numbers for granted, aims at constructing these objects from them. He then claims that, thanks to this proof, 
he would have succeeded ``in
defining the real number purely arithmetically or logically as a ratio of magnitudes that
are demonstrably there" (Frege 2013, p.~162$_2$). 

The adverb ``arithmetically'' is clearly used to emphasize 
that the envisaged definition would have been independent of both 
empirical considerations and geometry. 
In this sense, the definition would have surely been arithmetical, and our rendering of it will be as well. 
But there is another sense in which, despite his appealing to natural numbers, Frege did not
certainly want his definition to be arithmetical: both his criticisms to the alternative definitions 
depending on
an extension of the domain of rationals---including Cantor's (\S\S~II.68-85), Dedekind's (\S\S~II.138-147), and Weierstrass's
(\S\S~II.148-155)---and the very purpose of identifying real numbers with ratios of magnitudes make clear he 
wanted these numbers to be strictly independent of natural ones, to be 
properly \textit{Zahlen}, rather then \textit{Anzahlen}. By offering our definition, we will try, among
other things, to comply with this requirement.

In \S~6, we will
account for two distinct strategies
to get the required existence proof in our setting. One of them conforms to Frege's indications, while
the other might be considered more appropriate for ensuring logicality, since, \textit{pace} Frege, it does not require that the existence of continuously many objects be established. 
In \S~7, we will investigate whether the resulting definition does comply with logicality and non-arithmeticity---in the mentioned sense. In \S~8, we will provide some concluding remarks.

\section{Frege's Strategy\label{FS}}
Frege's strategy agrees with the ``application constraint'': the
requirement that a mathematical theory be shaped as to
immediately account for its applications\footnote{See Panza \& Sereni (2020) and Sereni (2019), which include a critical 
survey of the recent discussion on Frege's attitude toward applications of mathematical theories.}. This motivates his
suggestion to define real numbers as ratios of magnitudes, magnitudes as elements of distinct domains supposedly
including those of geometric, mechanic and empirical ones, and ratios on these domains as measures of the
relevant magnitudes. Insofar as it would be odd to require that the theory of real numbers involve these magnitudes as such,
together with their respective theories, this makes providing a structural definition of domains of magnitudes mandatory: a definition that merely fixes the conditions that a certain domain of
independent items has to meet in order to be recognized as a
domain of magnitudes. Frege himself clearly stresses this crucial point (Frege 1893-1903, \S II.161; Frege 2013, p.~158$_2$):
\begin{quotation}
There
are many different kinds of magnitudes: lengths, angles, periods of time, masses,
temperatures, etc., and it will scarcely be possible to say how objects that belong to
these kinds of magnitudes differ from other objects that do not belong to any kind
of magnitude. Moreover, little would be gained thereby; for we still lack any way of
recognizing which of these magnitudes belong to the same domain of magnitudes.

Instead of asking which properties an object must have in order to be a
magnitude, one needs to ask: how must a concept be constituted in order for
its extension to be a domain of magnitudes?
\end{quotation}

A natural way to
render the required structural definition would have provided definitional axioms, 
as usually done for groups or fields. An
informal conception of magnitudes recognizing the existence of ``lengths, angles, periods of time, masses,
temperatures, etc.''
might have suggested that there are non-isomorphic models satisfying these
axioms. Still, for Frege, magnitudes are  just those items that real numbers are
ratios of, and they all behave as lengths do, so that domains of magnitudes are all isomorphic to each other.
Had he defined them through appropriate axioms, these should have then been expected to be categorical, 
though algebraic in nature---as
it happens for the usual axioms for real numbers themselves, namely the
axioms of a totally ordered and Dedekind-complete field. Moreover, insofar as magnitudes
are required to add to each
other but not to multiply with each other (namely to admit only a single internal composition law), 
what he would have needed is a categorical axiomatization for
totally ordered, dense and Dedekind-complete (and, then, also Abelian and Archimedean)\footnote{%
See footnotes \ref{ft.GDCSAAA} and  \ref{ft.catMod} below.} groups.

Frege did not straightforwardly follow this route, however. Conforming
with a remark by
Gauss (1931, p. 635, also in \textit{Werke}, II, pp. 175-76; quoted in Frege, 1893-1903, \S~II.161) and putting it in his perspective,
he conceives of magnitudes as value-ranges of permutations, and so defines their domains not as domains
of items merely satisfying certain conditions, but rather as domains of extensions of appropriate first-level binary relations satisfying these conditions. This makes him able to appeal, along with his definitions, to structural properties of first-level binary relations, namely to the way
they compose and are inverted, as well as to there being the identity
relation among them. Accordingly, rather than listing a number of
axioms, and finally getting an implicit definition, Frege explicitly defines domains of magnitudes as
extensions of concepts under which extensions of certain first-level
relations fall. The difficulty he tackles is then that of looking for an explicit definition of the
concept of being one of these extensions (and, then, a magnitude), and falling under one of these concepts.

Since for Frege extensions are objects, this concept is first-level.
In order to define it, he appeals to a special function allowing him to reduce 
higher-level concepts to first-level ones, so as to work in a first-order fragment of his second-order theory. This is the
first-level two-arguments function $\xi \frown \zeta$, often too quickly identified with set-theoretic membership, whose definition is licensed by BLV. Once
BLV is omitted, this function can no more be defined, and the reduction to first-order is no more possible---unless by
a form of set theory. Hence, making 
Frege's definition consistent by eliminating BLV without falling into a set-theoretical setting requires replacing Frege's first-order definitions with higher-order ones. 
We will explain how this can be done by appropriately rephrasing Frege's definitions, and in clarifying
the logical nature of the (logical) system that is required for that. This is the purpose of
the next two \S\S.

\section{Frege's Definition of Domains of Magnitudes Rephrased}

\subsection{Eliminating Value-Ranges}

The omission of BLV is made possible by the elimination, from Frege's language, of terms for value-ranges. Insofar as
the presence of these terms in his definition of domains of magnitudes entirely depends on the function
$\xi \frown \zeta $, we have to make
its use pointless.

 To make a long story short, this function is such that for any objects $\Gamma $ and $\Delta $, if $\Gamma$ is the value-range $\overset{,}{\varepsilon }\Phi \left( \varepsilon \right) $ of a first-level one-argument function $\Phi \left( \xi \right) $\footnote{Our use
of Greek capital letters to denote objects and functions whatsoever corresponds to Frege's, in his 
``exposition'' of his formal language (Frege 1893-1903, Part I, \S~I.1-52).}, then $\Delta \frown \Gamma $ is $\Phi \left( \Delta \right) $, and if $\Gamma$ is not such a value-range---or, better, it is not a value-range at all, since, 
in Frege's formalism, any
value-range reduces to the value-range of a first-level one-argument function---, then $\Delta \frown \Gamma $ is the value-range of a first-level concept under 
which no object falls---for example that of $\lnot \left(\xi =\xi \right)$, which
we could denote by `$\oslash $', for short. In other terms, $\Delta \frown \Gamma $ is the value, for $\Delta $ as
argument, of the first-level one-argument function of which $\Gamma $ is the value-range, if $\Gamma$ is a value-range, and $\oslash $, if it isn't---whatever object $\Delta $ might be.

In a rich enough second-order predicate language including the operator `$\iota z\left[ z: \varphi \right] $' for definite descriptions, together with a symbol for value-ranges, the individual variables `$x$', `$y$' and `$z$', and the monadic predicate one `$F$', this stipulation could be rendered as follows:
\begin{equation*}
\forall x,y \left[ x\frown y = \iota z\left[z: \exists F\left( y=\overset{,}{\varepsilon }F\left(
\varepsilon \right) \wedge F\left( x\right) =z\right) \right]\right],
\end{equation*}%
provided that `$\iota z\left[z: \varphi \right]$' designates a well-defined object, namely $\oslash $, even
if there is no $z$ such that $\varphi $. If `$a$' and `$b$' are terms, this makes  `$a\frown b$' be a term in turn. 

This licenses using this term to denote the $\zeta$-argument of the same function $\xi \frown \zeta $.  Taking a new term `$c$' to denote the $\xi$-argument, one has the new term `$c\frown \left( a\frown b\right)$' such that
\begin{equation*}
c\frown \left( a\frown b\right)  = \iota z\left[ \exists G\left( b=\overset{,}{%
\alpha }\overset{,}{\varepsilon }G\left( \varepsilon ,\alpha \right) \wedge
G\left( c,a\right) =z\right) \right].
\end{equation*}%
It follows that `$c\frown \left( a\frown b\right)$' is a term that denotes the value for $c$ and $a$  as arguments
of the first-level two-argument function of which $b$ is the value-range, if $b$ is such a value-range\footnote{If $b$ is not
such a value-range, different cases are possible. It is not necessary to account for them, here. For a complete treatment, see
Panza (FC1).}.

Hence, if $\Phi\left( \xi \right)$ and 
$\Psi \left( \xi , \zeta \right)$ are a first-level one-argument and a first-level two-argument function, respectively, then
\[
a\frown \overset{,}{\varepsilon }\Phi
\left(\varepsilon \right) = \Phi \left(a \right)
\]
\[
c\frown \left(a \frown \overset{,}{\alpha } 
\overset{,}{\varepsilon }\Psi \left( \varepsilon, \alpha \right) \right) = \Psi \left(c, a \right).
\]

Suppose that `$\mathbf{P}_{b}$' and `$\mathbf{R}_{b}$' be respectively a monadic and a dyadic predicate\footnote{Here and in what follows, we take boldface capital latin letters as dummy letters for first-level properties and relations. The same letters in italics will, instead, be used for the corresponding variables.} appropriate 
for rendering, in an appropriate predicate language, 
two functions  $\Phi \left( \xi \right)$ and  $\Psi \left( \xi , \zeta \right)$ of which $b$ is the value-range. It follows that, in order to 
make the use of the function $\xi \frown \zeta $ pointless, and so eliminate value-ranges
while restating Frege's definition of domains of magnitudes
, it is enough to
replace each term of the form `$a\frown b$' with the formula `$\textbf P_{b} a $'
and each term of the form `$c\frown \left(a\frown b\right)$' with the formula `$c\mathbf R_{b} a$', and to
transform Frege's formal system accordingly\footnote{To be sure, this rendering of the function $\xi \frown \zeta $ in terms of predication is not fully faithful to Frege's original view: for Frege, both `$\Delta\frown\Gamma$' and `$\raisebox{-0,9mm}{{\mbox{---}}}\Delta\frown\Gamma$' (see footnote (\ref{ft.hor}), below) denote an object, while for us `$a\frown b$' is rather a formula. Nevertheless, whenever `$\Gamma$' denotes the extension of a concept $\Phi\left(\xi\right)$,  `$\Delta\frown \Gamma$' and `$\raisebox{-0,9mm}{{\mbox{---}}}\Delta\frown\Gamma$' are, for Frege, names of the True if and only if  `$\Phi(\Delta)$' itself is a name of the True. Insofar as only this case is relevant here, this rendering does not alter the aspects of Frege's
definition that are of interest here.}.
The system so obtained will be independent of BLV, and so will any definition stated in it.

\subsection{Working with Binary First-Order Relations}
Informally speaking, Frege conceived of a nonempty domain of magnitudes as a totally ordered, dense and
Dedekind-complete additive group of permutations.
In light of his rejection of implicit definitions, defining such a group required to
explicitly defining a particular function to play the role of its (additive) law
of composition, which required, in turn, to have objects available, endowed with an internal structure making
such a definition possible. To this purpose, he made the simplest choice possible: 
he took those objects to
be extensions of functional first-level binary relations, and assigned the role of this law
to the composition of the
corresponding relations. This obviously resulted in taking the extension of the identity relation as the
neutral element of the group, and
the extensions of the inverse relations as its inverse elements.

This choice is easily rendered in a predicate setting lacking extensions. We merely
have to fix the conditions under which a first-level binary relation is 
functional and results either from the inversion of another such relation,  or from the composition of two other such relations.
Taking `$R$' and `$S$' to range over first-level binary relations, this can be formally done through the following explicit definitions:
\[
\begin{array}{rl}
\left[\text{Functionality}\right] & \forall R \left( \mathscr{I} R \Leftrightarrow \forall x,y (xRy \Rightarrow
\forall z (xRz\Rightarrow y=z))\right),\medskip \\
\left[\text{Inversion}\right] & \forall R \forall x,y \left( xR^-y \Leftrightarrow yRx \right), \medskip \\
\left[\text{Composition}\right] & \forall R,S \forall x,y \left(x[R \sqcup S]y \Leftrightarrow \exists z \left(xRz \wedge zSy\right)
\right).
\end{array}
 \]

These definitions define three third-order constants, respectively: the monadic predicate 
constant `$\mathscr{I}$', designating a property of first-level binary
relations; the monadic functional constant `$\sim ^-$', designating a
one-argument function from and to first-level binary relations; the dyadic functional constant `$\sim \sqcup \sim$', designating a
two-argument function from and to first-level binary relations, again. Hence, in
order to be licensed, they require appropriate instances of predicative
comprehension. The first 
requires the following instance of third-order
predicative comprehension without parameters: 
\begin{equation}
\exists \mathscr X
\forall R \left(\mathscr XR\Leftrightarrow \forall x,y(xRy \Rightarrow \forall z
(xRz\Rightarrow y=z))\right).\tag{Functionality-CA}
\end{equation}
The other two
require the following
instances of second-order dyadic predicative comprehension with parameters respectively: 
\begin{equation}
\forall R \exists S \forall x,y \left(xSy\Leftrightarrow yRx\right),\tag{Inversion-CA}
\end{equation}
\begin{equation}
\forall R,S
\exists T \forall x,y \left( xTy \Leftrightarrow \exists z (xRz\wedge zSy))\right),\tag{Composition-CA}
\end{equation}
where `$T$' ranges over first-level binary relations, too.

One might replace, however, these explicit definitions with 
mere (metalinguistic)
typographic stipulations:
\[
\begin{array}{rl}
\left(\text{Functionality$'$}\right) & \mathscr I(R) := \forall x,y (xRy \Rightarrow \forall z
(xRz\Rightarrow y=z)),\medskip \\
\left(\text{Inversion$'$}\right) & R^-\left(xy\right) :=yRx, \medskip \\
\left(\text{Composition$'$}\right) & R\sqcup S\left(xy\right) := \exists z (xRz\wedge Szy).
\end{array}
 \]
Any instance of the left-hand side of these
stipulations is intended to be 
a mere abbreviation of the corresponding
instance of the right-hand side. For example, while
`$\mathscr{I}R $' in (Functionality) is an atomic third-order
(open) formula, 
`$\mathscr I(R)$' in (Functionality$'$) is an atomic symbol abbreviating the second-order (open) 
formula $\text{`}\forall x,y (xRy\Rightarrow \forall z
(xRz\Rightarrow y=z)) \text{'}.$ And analogously for  `$xR^-y$' and `$R^-\left(xy\right)$' 
in (Inversion) and (Inversion$'$), respectively, and for `$x\left[R\sqcup S\right]y$' and `$R\sqcup S\left(xy\right)$' in (Composition) and (Composition$'$), respectively. Adopting
these stipulations requires neither any instance of comprehension, nor any extension of the usual second-order language.

We will see in what follows whether these stipulations are enough for our purpose,
or the corresponding explicit
definitions are needed, and the instances of comprehension they require\footnote{In Frege's original setting things would not be so simple. Consider only the example of (Functionality). In this setting, the role of this definition is played by the definition of the first-level concept $\textbf{I}\xi$ (Frege, 1893-1903, \S~I.37). By adapting Frege's notation to our modern one, the definition might be stated as follows:
\begin{equation*}
 \left[\forall x,y  \left[\raisebox{-0,9mm}{{\mbox{---}}}\left( x \frown \left( y \frown a\right)\right) \Rightarrow 
\forall z\left[\raisebox{-0,9mm}{{\mbox{---}}}\left( x \frown \left( z \frown a\right)\right) \Rightarrow y=z  \right]\right]\right]=%
\textbf{I}a.
\end{equation*}
where `$a$' is a term used as a parameter, and $\raisebox{-0,9mm}{{\mbox{---}}}\xi$ is the horizontal concept (\textit{ibidem}, \S~I.8), which is such that $\raisebox{-0,9mm}{{\mbox{---}}}\Gamma$
is the True if $\Gamma$ is also the True, and the False otherwise. 
It follows that $\textbf{I}a$ is the same object as 
$\forall x,y  \left[\raisebox{-0,9mm}{{\mbox{---}}}\left( x \frown \left( y \frown a\right)\right) \Rightarrow 
\forall z\left[\raisebox{-0,9mm}{{\mbox{---}}}\left( x \frown \left( z \frown a\right)\right) \Rightarrow y=z\right]\right]$,
which is a truth-value.  If $a $ is not a value-range of a first-level binary relation, $\raisebox{-0,9mm}{{\mbox{---}}}\left(b \frown \left( c \frown a\right)\right)$ is the False for whatever pair of objects $b$ and $c$, and 
$\textbf{I}a $ is then the True, which makes any object other than a value-range of a first-level binary relation fall under the concept $\textbf{I}\xi $. 
If $a $ is
the value-range of a first-level binary relation $\Phi \left( \xi ,\zeta \right) $, $a$ falls under the concept $\textbf{I}\xi $ if and only if
either  $\Phi \left( \xi ,\zeta \right) $ is empty, or, for any $x$, there is at most one  $y$ such that $\Phi \left(x ,y \right) $ is the True. 
Clearly, there is no way to regard this definition as a mere typographic stipulation. 
It rather defines a total first-level concept by introducing a functional constant to designate it.
Among many others, there are two relevant differences with our case: \textit{i})~Frege's definition applies in general, whereas 
both (Functionality)  and (Functionality$'$) only
apply to first-order binary relations; \textit{ii})~Differently from (Functionality$'$), Frege's definition is licensed only via 
a stipulation analogous to second-order comprehension.
\textit{Mutatis mutandis}, this also applies to (Inversion) and (Composition), and to any other particular definition entering his definition of domains of magnitudes.\label{ft.hor}}.

\subsection{Domains of Classes}
For Frege, a domain of magnitudes is the domain of a ``positive class'', which is in turn 
a ``positival class'' of an appropriate sort. 
In his setting, a class is the extension of a first-level concept (Frege, 1893-1903, \S~II.16),
and the objects falling under this concept are said to belong to the class.
Positival and positive classes are, in particular, extensions of concepts under
which (only) extensions of first-level
binary relations fall. Defining them amounts to 
fixing the conditions 
that a concept is to meet for the objects falling under it to
be just these extensions.
To do this, Frege appeals to their ``domains''. He has, then, to firstly define, in general,
domains of classes (\textit{ibidem}, \S~II.173). The definition applies to any class, but we only need to consider its application
to the case of the domain of a class of extensions of first-level binary relations. 

This is the extension of a concept under which fall: the extensions in the class; 
the extensions of the inverses of the relations whose extensions are in the class; and the extensions of the relations composed by
each of the relations whose extensions are in the class 
and their inverses---which in case these relations are
functional, as required for both positival and positive
classes, all coincide with the extension of the identity relation.
In our setting, we can, then, rephrase, Frege's definition of the domain of a class of extensions of first-level 
binary relations as follows:
\begin{equation}
\forall \mathscr X \forall R \left( \eth \mathscr X R \Leftrightarrow
\left\{ 
\begin{array}{l}
\mathscr X R \vee \\ 
\exists S\left[ \mathscr X S\wedge \forall x,y\left[ 
\begin{array}{l}
\left[xRy\Leftrightarrow S^-\left(xy\right)\right] \vee \\ 
\left[xRy\Leftrightarrow S\sqcup S^-\left(xy\right)\right]%
\end{array}%
\right] \right]%
\end{array}%
\right\} \right)\label{3.4}
\end{equation}
where `$\mathscr X $' is a third-order monadic variable, and `$\eth $' a functional operator applied to it.
This definition makes clear that, when applied to whatever
(second-level) property $\mathscr Q$ of first-level binary relations\footnote{Here
and in what follows we use  `$\mathscr Q$' as a dummy letter for second-level properties. Later we will also use  `$\mathscr A$', `$\mathscr E$', `$\mathscr H$',
`$\mathscr L$', `$\mathscr M$'  and `$\mathscr P$' for the same purpose.},
$\eth $ gives another property $\eth \mathscr Q$ of these same relations.

To license this definition, we need to ensure the existence and uniqueness
of a second-level property providing a putative value for $\eth \mathscr X$ 
under the existence of a second-level property providing a value for 
$\mathscr X$, and this requires, in turn, third-order comprehension with parameters.
But suppose we wanted to define a certain (third-level) property $\mathcal Q$\footnote{Here
we use  `$\mathcal Q$' as a dummy letter for third-level properties.} that a
class of first-level binary relations should have in order to be positival, which
is required to render Frege's definition of positival classes. If, in defining it, we had to appeal to the
domains of the classes that could have it, as is also required to render Frege's
definition, we should have
recourse to a definition like this: 
\begin{equation*}
\forall \mathscr X \left[ \mathcal{Q }\mathscr X \Leftrightarrow \phi
\left(\eth \mathscr X \right) \right]
\end{equation*}
where `$\phi \left(\eth \mathscr X \right)$' stands for an
appropriate formula involving the predicate `$\eth \mathscr X $'. Hence, insofar as, in our
rendering of Frege's definition of positival and positive classes and domains of magnitudes, this predicate would
only appear in instances of formulas of the form `$\eth \mathscr X R$', we
can replace (\ref{3.4}) with the following abbreviation stipulation
\begin{equation}
\eth (\mathscr X)(R) := \left\{ 
\begin{array}{l}
\mathscr X R \vee \\ 
\exists S\left[ \mathscr X S\wedge \forall x,y\left[ 
\begin{array}{l}
\left[xRy\Leftrightarrow S^-\left(xy\right)\right] \vee \\ 
\left[xRy\Leftrightarrow S\sqcup S^-\left(xy\right)\right]%
\end{array}%
\right] \right]%
\end{array}%
\right\}, \tag{\ref{3.4}$'$}
\end{equation}
then use appropriate instances of `$\eth (\mathscr X)(R)$' instead of the corresponding instances of `$\eth %
\mathscr X R$'. As a matter of fact, this stipulation is all we need for our present purpose, and it must be supplied 
by no sort of comprehension, since it
introduces no new predicate, but merely lets each instance of its
left-hand side be an abbreviation of the corresponding
instance of the right-hand side. For short, read both `$\eth X R$' and `$\eth (\mathscr X)(R)$' as `$R$ belongs to the domain of the class of first-order binary relations that have $\mathscr X$'.

\subsection{Positival Classes\label{Sect.PC}}
We can now consider Frege's definition of positival classes. If we had to render it 
through a(n explicit) definition, we should define
a fourth-order monadic predicate constant designating
a third-level property. This would require to
quantify over second-level properties, and, then, to appeal to fourth-order
comprehension. But, once again, we are not forced to do it.
As above, we might 
recur to an abbreviation stipulation by so avoiding any sort of comprehension.

In agreement with Frege's definition, the extension of a first-level binary relation $\mathbf R$
belongs to a positival class if (and only if): both $\mathbf R$ and its inverse are
functional; the extension of $\mathbf R \sqcup \mathbf R^{-}$, i.e. the
identity relation, does not belong to
the class; and for any first-level binary relation $\textbf S$, if its extension belongs to
the class, then: the class of the objects that bear $\mathbf R$ to some other object
coincides with the class of the objects to which some object bears $\textbf S$\footnote{This makes the relations whose 
extensions belong to the class permutations on a subjacent
first-order domain.}; the extension of $\mathbf R\sqcup \textbf S$ belongs to the class; both
the extension of $\mathbf R^{-}\sqcup \textbf S$ and that of $\mathbf R\sqcup \textbf S^{-}$ belong to the
domain of the class. In our setting, this can be rendered either this way
\begin{equation}
\forall \mathscr X \left[\mathcal{L} \mathscr X \Leftrightarrow \forall R\left[ \mathscr X R\Rightarrow %
\left[ 
\begin{array}{l}
\forall S\left[ \mathscr XS\Rightarrow \left[ 
\begin{array}{l}
\forall x\left[ \exists y\left( xRy\right) \Leftrightarrow \exists z\left(
zSx\right) \right] \wedge \medskip \\ 
\mathscr XR\sqcup S\wedge \medskip \\ 
 \eth \mathscr X  R\sqcup S^{-} \wedge \medskip \\ 
\eth \mathscr X R^{-}\sqcup S%
\end{array}%
\right] \right] \wedge \medskip \\ 
\mathscr{I}R \wedge \mathscr{I}R^{-} \wedge
\lnot \mathscr XR\sqcup R^{-}%
\end{array}%
\right] \right]\right],\label{3.2}
\end{equation}
or this way:
\begin{equation}
\mathcal{L} (\mathscr X) :=\forall R\left[ \mathscr X R\Rightarrow %
\left[ 
\begin{array}{l}
\forall S\left[ \mathscr XS\Rightarrow \left[ 
\begin{array}{l}
\forall x\left[ \exists y\left( xRy\right) \Leftrightarrow \exists z\left(
zSx\right) \right] \wedge \medskip \\ 
\mathscr XR\sqcup S\wedge \medskip \\ 
\eth (\mathscr X) \left( R\sqcup S^{-}\right) \wedge \medskip \\ 
\eth (\mathscr X) \left( R^{-}\sqcup S\right)%
\end{array}%
\right] \right] \wedge \medskip \\ 
\mathscr{I}(R) \wedge \mathscr{I}(R^{-}) \wedge
\lnot \mathscr XR\sqcup R^{-}%
\end{array}%
\right] \right] \tag{\ref{3.2}$'$},
\end{equation}
where both `$\mathcal{L} \mathscr X$' and `$\mathcal{L}\left( \mathscr X\right)$' are short for `$\mathscr X$ is a positival class' or, more
precisely, `the first-level binary relations having $\mathscr X$ form a positival class'. 

In (\ref{3.2}), `$\mathcal{L}$' is a fourth-order  predicate constant and `$\mathcal{L} \mathscr X$' an atomic (open) formula.
This is, then, an explicit definition, which is to be 
licensed by an appropriate form of fourth-order comprehension. In (3.2$'$), `$\mathcal{L} \left(\mathscr X\right)$' is, instead,
an abbreviated (open)
formula, and `$\mathcal{L}$'  is merely a symbol occurring  in it. Hence (3.2$'$) neither requires a fourth-order language nor is 
to be licensed by any form of fourth-order 
comprehension.

This does not mean that no form of comprehension is required to license it.
Since, in its right-hand side, the signs `$^{-}$' and `$\sqcup $' do not merely occur as
parts of the abbreviated formulas `$R^{-}\left( xy\right) $' and `$R\sqcup
S\left( xy\right) $' introduced by  
(Inversion$'$) and (Composition$'$), but as functional signs
allowing to construe the predicate variables `$R^{-}$' and `$R\sqcup S$'. Their use in (3.2$'$) is, then, to be 
licensed by the explicit
definitions (Inversion) and (Composition), which respectively require, in turn,
(Inversion-CA) and 
(Composition-CA), or, more generally, the following second-order predicative comprehension axiom schema with parameters:
\begin{equation}
\forall R\ldots T
\exists U \forall x,y \left[ xUy \Leftrightarrow \phi_{\Delta_0^1}\left(R \ldots T\right)\right],\tag{PCA$^{2}_{\Delta^{1}_{0}}$}
\end{equation} 
where `$U$', `$R$', `$S$', and `$T$' range over first-level binary relations,
and  `$ \phi_{\Delta_0^1}\left(R \ldots T\right)$' stands for any second-order formula containing the parameters `$R$', \ldots `$T$', but no higher-order quantifiers.
 
Before going ahead with the definition of positive classes, some remarks are in order about the informal import of the
conditions characterizing a positival class. They apply, \textit{mutatis mutandis}, both to (\ref{3.2}) and to  (\ref{3.2}$'$), but, 
for short and simplicity, we only make them about the latter.

Let $\mathscr L$ be a second-level property. Requiring that 
\[
\forall R\left[ \mathscr L R\Rightarrow \left( \mathscr{I}(R) \wedge \mathscr{I}
(R^{-}) \right) \right]
\]
amounts to requiring that both a binary relation
that has $\mathscr L$ and its inverse are functional. If this condition obtains, requiring that  
\[
\forall R\left[\mathscr L R\Rightarrow \lnot \mathscr L R\sqcup
R^{-} \right]
\] 
and that 
\[
\forall R,S \forall x \left[ \left(\mathscr L R \wedge \mathscr L S\right) \Rightarrow \exists
y\left[\left( xRy\right) \Leftrightarrow \exists z\left( zS%
x\right) \right] \right]
\]
respectively amount to requiring that the identity relation has
not $\mathscr L$, and that all the relations having $\mathscr L$ are permutations%
\footnote{%
One should better say `correspond to permutations', since, strictly 
speaking, permutations are functions, not relations. Let us adopt,
however, a more straightforward, though abusive, language, for short.} on a subjacent unspecified set.
Hence, only permutations but the identity one, have $\mathscr L$. 
Thus, $\sqcup $ is an
associative law of composition without neutral element on the relations having $\mathscr L$. 
Again, if
all the above conditions obtain, requiring that
\[
\forall R,S \left[\left(\mathscr L R \wedge \mathscr L S\right) \Rightarrow 
\mathscr LR\sqcup S \right]
\]
amounts to requiring that the family of permutations having $\mathscr L$ is closed
under $\sqcup $. This makes: the inverse of any such permutation not have $\mathscr L$---since, if it did,  the identity permutation would also have it; the family of permutations that 
satisfy the open formula `$\eth  (\mathscr L) (R) $' be also 
closed under composition of the inverses of those having $\mathscr L$---since, for whatever permutations $\mathbf{R}$ 
and $\mathbf{S}$ that
have  $\mathscr L$,
$\mathbf{R}^{-}\sqcup \mathbf{S}^{-}$ is the same
permutation as $\left( \mathbf{S}\sqcup \mathbf{R}\right) ^{-}$. All this
is still not enough to ensure that the family of permutations that 
satisfy the open formula `$\eth  (\mathscr L) (R) $', if any, is closed under $\sqcup $,
and forms, then, a(n additive) group of permutations. Also requiring that
\[
\forall R,S \left[\left(\mathscr L R \wedge \mathscr L S\right) \Rightarrow \left[\eth (\mathscr L) \left(R \sqcup S^{-}\right) 
\wedge \eth (\mathscr L) \left( R^{-}\sqcup S\right)\right]\right]
\]  
just amounts to requiring it.
If $\mathscr L$ is a second-level (monadic) property such that $\mathcal{L}\left(\mathscr L \right)$, the
first-level binary relations satisfying the open formula `$\eth  (\mathscr L) (R) $', if any, form, then, a(n additive) group of
permutations, whose 
internal law of composition is $\sqcup $,  whose neutral element is the identity permutation, and whose  inverse function is 
$R \longmapsto R^{-}$.

This group is not necessarily Abelian, for $\sqcup $ is not commutative on permutations.
But it is endowed with a
total and right-invariant order defined in terms of the composition
operation. Since, if $\mathbf{H}$ and $\mathbf{K}$ are two permutations
whatsoever that satisfy `$\eth \left( \mathscr L\right) \left( R\right) $',
requiring that $\mathbf{H}\sqcup \mathbf{K}^{-}$ have the property $\mathscr %
L$ is equivalent to requiring that $\mathbf{K}$ and $\mathbf{H}$ bear a 
right-invariant strict-order relation, let as say $\sqsubset_{\mathscr L}$, on these permutations\footnote{%
The proof is simple. As it has been required that $\lnot \mathscr L\mathbf{H}\sqcup \mathbf{H}%
^{-}$, we immediately have that  $\lnot \mathbf{H} \sqsubset_{\mathscr L} \mathbf{H}$. As 
$\mathbf{K}\sqcup \mathbf{H}^{-}$ is the same permutation as $\left( \mathbf{H}\sqcup \mathbf{K}%
^{-}\right) ^{-}$, we have that 
$\mathscr L\mathbf{H}\sqcup \mathbf{K}^{-}\Rightarrow \lnot %
\mathscr L\mathbf{K}\sqcup \mathbf{H}^{-}$, i.e.
$ \mathbf{K} \sqsubset_{\mathscr L} \mathbf{H} \Rightarrow \lnot \mathbf{H} \sqsubset_{\mathscr L}  \mathbf{K}$.
Again, if $\mathbf{J}$ is, also, a permutation
that satisfies `$\eth \left( \mathscr L\right) \left( R\right) $',
then $\mathbf{J}\sqcup \mathbf{K%
}^{-}$ is the same permutation as $\left( \mathbf{J}\sqcup \mathbf{H}%
^{-}\right) \sqcup \left( \mathbf{H}\sqcup \mathbf{K}^{-}\right) $, and so we
have that $\left( \mathscr L\mathbf{H}\sqcup \mathbf{K}%
^{-}\wedge \mathscr L\mathbf{J}\sqcup \mathbf{H}^{-}\right) \Rightarrow %
\mathscr L\mathbf{J}\sqcup \mathbf{K}^{-}$, i.e.
$\left(\mathbf{K} \sqsubset_{\mathscr L} \mathbf{H} \wedge \mathbf{H} \sqsubset_{\mathscr L} \mathbf{J}\right)   \Rightarrow
\mathbf{K} \sqsubset_{\mathscr L} \mathbf{J}$.
Finally, as  $\left( \mathbf{H}\sqcup \mathbf{J}\right) \sqcup
\left( \mathbf{K}\sqcup \mathbf{J}\right) ^{-}$ is the same permutation
as $\mathbf{H}\sqcup \mathbf{K}^{-}$, we have that  $%
\mathscr L\mathbf{H}\sqcup \mathbf{K}^{-}\Rightarrow \mathscr L\left( 
\mathbf{H}\sqcup \mathbf{J}\right) \sqcup \left( \mathbf{K}\sqcup \mathbf{J}%
\right) ^{-}$, i.e. $\mathbf{K} \sqsubset_{\mathscr L} \mathbf{H}
\Rightarrow \mathbf{K} \sqcup \mathbf{J} \sqsubset_{\mathscr L} \mathbf{H} \sqcup \mathbf{J}$.}. Hence, 
if this relation
is conceived of as the smaller-than relation (that is, `$\mathscr L\mathbf{H}%
\sqcup \mathbf{K}^{-}$' or `$\mathbf{K} \sqsubset_{\mathscr L} \mathbf{H}$' are read as
 `$\mathbf{K}$ is smaller than $\mathbf{H}$%
'), then we can take the collection of the permutations that have $\mathscr L
$, if any, as the positive semi-group of the group of permutations formed by the permutations
that satisfy `$\eth (\mathscr L)(R)$'\footnote{%
In commenting his definition of positival classes, Frege (1893-1903, \S %
~II.175; Frege 2013, pp.\ 171$_{2}$-72$_{2}$) claims to have \textquotedblleft tried [%
\textit{bem\"{u}ht}]\textquotedblright\ to include in it only
\textquotedblleft necessary [\textit{nothwendigen}]\textquotedblright\ and
\textquotedblleft mutually independent [\textit{einander unabh\"{a}ngig}%
]\textquotedblright\ conditions, though taking as unprovable his having
succeeded in this. In a note added at the end of his book (\textit{ibidem} vol.\ 2, p.
243, Frege 2013, p.~243$_{2}$), explicitly referred to this comment, he corrects
himself by observing that a proof could have been possible by means of
counterexamples, though taking it to be \textquotedblleft doubtful [\textit{%
bezweifeln}]\textquotedblright\ that these counterexamples could be given
in his formal setting. Dummett (1991, p. 288) suggests that his
doubt concerned the independence of the condition we
expressed by `$\forall R,S \left[\left(\mathscr L R \wedge \mathscr L S\right) \Rightarrow
\eth (\mathscr L) \left( R^{-}\sqcup S\right)\right]$' from the other ones characterizing 
a positival class, by observing that, in his developments concerning domains of magnitudes,
Frege appeals to this condition as
late as possible (namely only in 
\S ~II.218), after making explicit (\S~II.217) the \textquotedblleft
indispensability\textquotedblright\ of this condition for the purpose for
which it is used, which, in our setting, corresponds to prove that if $\mathbf{H}$ and $\mathbf{K}$ belong
to a positive class and $\mathbf{H}$ is smaller than $\mathbf{K}$ over the positive semigroup
involved in this class, then $\mathbf{K}^{-}$ is smaller than $\mathbf{H}^{-}$ over the corresponding group.
Adeleke, Dummett and Neumann (1987, th. 2.1)  
have finally proved that this condition is actually
independent of the others.
When transposed in our setting, the proof goes along the following lines. Let $\mathscr L^{\star}$ be a property
satisfying `$\mathcal L \left(\mathscr X \right)$' except for the condition at issue, 
$\mathcal{G^{\star}}$ be the structure
formed by the permutations that satisfy `$\eth \left( \mathscr %
L^{\star}\right) \left( R\right)$, and 
$\mathbf{H}$ and $\mathbf{K}$ two binary first-level relations having 
$\mathscr L^{\star}$. Insofar as $\left( \mathbf{K}^{-}\sqcup \mathbf{H}\right)
^{-}$ is the same permutation as $\mathbf{H}^{-}\sqcup \mathbf{K}$,
not ensuring that $\mathbf{K}^{-}\sqcup \mathbf{H}$ satisfies `$\eth \left( %
\mathscr L^{\star}\right) \left( R\right) $' is the same as not ensuring that the
disjunction 
\begin{equation*}
`\mathscr L\mathbf{K}^{-}\sqcup \mathbf{H}\vee \forall x,y\left[
x\mathbf{H}^{-}y\Leftrightarrow x\mathbf{K}^{-}y\right] \vee \mathscr L%
\mathbf{H}^{-}\sqcup \mathbf{K}\text{'} \qquad \text{i.e.} \qquad `\mathbf{H}^{-}\sqsubset_{\mathscr L} \mathbf{K}^{-}\vee \mathbf{H}^{-}=_{\mathscr L}\mathbf{K}^{-} \vee
\mathbf{K}^{-}\sqsubset_{\mathscr L} \mathbf{H}^{-}\text{'}
\end{equation*}%
holds, namely that $\mathbf{H}^{-}$ and $\mathbf{K}^{-}$ are comparable
according to the order over $ \mathcal G^{\star}$. It would follow that, besides of not being a group, this last structure
is not endowed with a total order, but only with a
partial one. It can be proved (\textit{ibidem}, Lemma 1.2) 
that this partial order is an
\textquotedblleft upper semilinear order\textquotedblright---that is, a strict
partial order ``such that the elements greater than any given one are
comparable, and that, for any two incomparable elements, there is a third
greater than both of them\textquotedblright, or, more simply, a strict partial
order that \textquotedblleft may branch downwards, but cannot branch
upwards\textquotedblright\ (Dummett 1991, p. 288). 
But $\mathcal{G}^{\star}$ is a sub-structure of 
the group $\mathcal{G}$ formed by the permutations that satisfy `$\eth \left( \mathscr %
L\right) \left( R\right)$, where $\mathscr L$ satisfies `$\mathcal L \left(\mathscr X \right)$' as a whole. 
Hence, the condition at issue follows from the others if
and only if $\mathcal{G}^{\star}$ can be extended in no group
other than $\mathcal{G}$. To prove the independence of this condition 
it is, then, enough to prove that there is a group including $\mathcal{G}^{\star}$ other than $\mathcal{G}$.
By Cayley's theorem, any group is isomorphic to a group of permutations. It is, then, enough to prove that
there is a partially ordered group
whatsoever not isomorphic to $\mathcal{G}$ (that is, not totally ordered) that includes
a sub-structure isomorphic to $\mathcal{G}^{\star}$.
This is just what Adeleke, Dummett and Neumann do. \label{ADMQF}}.

\subsection{Positive Classes and Domains of Magnitudes\label{Sect.PCDM}}
Informally speaking, a nonempty positive
class is a positival class whose domain is a totally-ordered, dense and
Dedekind-complete
group of
permutations, which is, by consequence,  also Archimedean and Abelian\footnote{%
That a totally-ordered,
dense and Dedekind-complete group (of permutations) is also Archimedean and
Abelian is in fact proved by Frege himself. He proves that a Dedekind-complete positival class is Archimedean
(Frege 1893-1903, th. 635, \S\ II.~214), and that the domain of a positive class is Abelian
(\textit{ibidem}, th. 689, \S\ II.~244). This is the last theorem he proves.
Insofar as the proof of the former theorem does not appeal to the condition
considered in footnote (\ref{ADMQF}) above, Adeleke, Dummett and Neumann (1987, p.~516)
restate these theorems as follows: a Dedekind-complete upper semilinear
order is Archimedean---which, of course makes it also a Dedekind-complete
total order; if the order of a group is dense, Archimedean and total, then
the group is Abelian.\label{ft.GDCSAAA}}. By having a
strict order available, the density condition can be stated easily. For stating the  Dedekind-completeness one,
further means are required. 

To this purpose, Frege defines
the upper rims over a
positival class. Let $\mathscr L$ be such that $\mathcal L \left(\mathscr L \right)$, and $\mathscr A$ a sub-property
of it. In our setting, an upper rim $\mathbf{U}$ of the collection of
permutations that have $\mathscr A$ over the collection of those that have $%
\mathscr L$ is a relation having $\mathscr L$, such that any other relation
that has $\mathscr L$ and is smaller than $\mathbf{U}$ over the former
collection has $\mathscr A$. To define it, Frege begins with a general
definition, then
applies it to positival classes. In the general case, both the 
informal notion of an upper rim and the subsequent one of an upper limit
become nonsensical. Formally speaking, this is immaterial, however, since the following definition of positive classes
excludes that the deviant cases obtain in the case of such a class\footnote{See footnote \ref{ft.URUL} below.}.

Once again, there are two ways to render the general definition: either as
\begin{equation}
\forall \mathscr X, \mathscr Y \forall R \left[\mathscr X \Finv_{R} \mathscr Y
\Leftrightarrow 
\forall S\left[ \left( \mathscr XS \wedge \mathscr XR\sqcup
S^{-}\right) \Rightarrow \mathscr YS\right]\right],\label{3.6}
\end{equation}
or as
\begin{equation}
\left[ \left( \mathscr X\right) \Finv \,\left( \mathscr Y\right) \right]
\left(R\right) :=\forall S\left[ \left( \mathscr XS\wedge \mathscr X(R\sqcup
S^{-}\right)) \Rightarrow \mathscr YS\right]. \tag{\ref{3.6}$'$}
\end{equation}
 
The upper limit of a sub-class of a positival class is the greatest of all the upper rims of the former
over the latter, if
there is one. Anew, Frege's definition can be rendered in two ways: either as
\begin{equation}
\forall \mathscr X, \mathscr Y \forall R \left[\mathscr X \text{\l}_{R} \mathscr Y
\Leftrightarrow
\left\{ 
\mathcal{L} \mathscr X \wedge \mathscr X R \wedge \mathscr X \Finv_{R} \mathscr Y \wedge
\neg\exists S\left[ 
\mathscr XS \wedge \mathscr XS \sqcup R^{-} \wedge 
\mathscr X \Finv_{S} \mathscr Y
\right]
\right\}\right],\label{3.7}
\end{equation}
or as:
\begin{equation}
[\left(\mathscr X\right) \text{\l } \left(\mathscr Y\right)]\left(R\right):=
\left\{ 
\mathcal{L} \left(\mathscr X\right) \wedge \mathscr X R \wedge [\left(
\mathscr X\right) \Finv \left(\mathscr Y\right)]\left(R\right) \wedge
\neg\exists S\left[ 
\begin{array}{l}
\mathscr XS \wedge \mathscr X(S \sqcup R^-) \wedge 
\medskip \\ 
\left[\left(\mathscr X\right) \Finv \left( \mathscr Y\right)\right]%
\left(S\right)
\end{array}
\right]
\right\}\tag{\ref{3.7}$'$}.
\end{equation}

For short, read both `$\mathscr X \Finv_{R} \mathscr Y$'
and `$\left[\left(\mathscr X\right) \Finv \, \left(\mathscr %
Y\right)\right] \left(R\right)$' as
`$R$ is an upper rim of $\mathscr Y$ over~$\mathscr X$', and both `$\mathscr X \text{\l}_{R} \mathscr Y$'
and `$[\left(\mathscr X\right) \text{\l } \left(\mathscr Y\right)]\left(R\right)$' as
`$R$ is the upper limit of $\mathscr Y$ over~$\mathscr X$'.

Let $\mathscr P$ be a third-level monadic property of first-order binary relations. 
Informally speaking, the relations that have it form a positive class if they form a positival one, and are such that: for any relation $R$ which has $\mathscr P$, there is another relation $S$ smaller than it over $\mathscr P$ (density); any proper subclass $ \mathscr Y $ of $ \mathscr P $ which has an upper rim over $ \mathscr P $ also has  an upper limit over $ \mathscr P $ (Dedekind-completeness). These conditions can be rendered in two ways: either as
\begin{equation}
\forall \mathscr X \left[
\mathcal{P} \mathscr X \Leftrightarrow \left\{ 
\begin{array}{l}
\mathcal{L} \mathscr X \wedge \medskip \\ 
\forall R\left[ \mathscr XR \Rightarrow \exists S\left[ \mathscr XS \wedge %
\mathscr X R\sqcup S^-\right] \right] \wedge \medskip \\ 
\forall \mathscr Y \left[ \left[ 
\begin{array}{l}
\exists R\left[ \mathscr X \Finv_{R} \mathscr Y \wedge \mathscr X R \right] \wedge \medskip \\ 
\exists S\left[ \mathscr X S \wedge \lnot \mathscr Y S \right]%
\end{array}
\right] \Rightarrow \exists T\left[ \mathscr X \text{\l}_{T}\mathscr Y\right] \right]%
\end{array}%
\right\}\right],  \label{3.8}
\end{equation}
or as:
\begin{equation}
\mathcal{P} \left(\mathscr X\right) := \left\{ 
\begin{array}{l}
\mathcal{L} \left(\mathscr X\right) \wedge \medskip \\ 
\forall R\left[ \mathscr XR \Rightarrow \exists S\left[ \mathscr XS \wedge %
\mathscr X R\sqcup S^-\right] \right] \wedge \medskip \\ 
\forall \mathscr Y \left[ \left[ 
\begin{array}{l}
\exists R\left( \left[\left(\mathscr X\right) \Finv \left(\mathscr Y\right)%
\right]\left(R\right) \wedge \mathscr X R \right) \wedge \medskip \\ 
\exists S\left( \mathscr X S \wedge \lnot \mathscr Y S \right)%
\end{array}
\right] \Rightarrow \exists T\left[ \left[\left(\mathscr X\right) \text{%
\textrm{\l }}\left(\mathscr Y\right)\right] \left(T\right) \right] \right]%
\end{array}%
\right\}.  \tag{\ref{3.8}$'$}
\end{equation}

For short, read both `$\mathcal{P} \mathscr X$' and `$\mathcal{P}\left( \mathscr X\right)$' as `$\mathscr X$ is a positive class' or, more
precisely, `the first-level binary relations having $\mathscr X$ form a positive class'.

While (\ref{3.6}), (\ref{3.7}) and (\ref{3.8}) are explicit definitions, and have to be 
licensed by some form of fourth-order comprehension, (\ref{3.6}$'$), (\ref{3.7}$'$) and (\ref{3.8}$'$) are abbreviation stipulations, 
and require no form of comprehension stronger than
(PCA$^{2}_{\Delta^{1}_{0}}$)\footnote{\label{ft.URUL}Comparing (\ref{3.8}) and (\ref{3.8}$'$), on the one side, with (\ref{3.6}-\ref{3.7}) and 
(\ref{3.6}$'$-\ref{3.7}$'$), on the other, allows us to see why the deviant cases pertaining to the definition of an upper rim
and an upper limit become immaterial by passing to the definition of a positive class. For short and simplicity, we only consider 
(\ref{3.6}$'$-\ref{3.8}$'$).
The right-hand side of (\ref{3.6}$'$) fails, as such, in rendering the informal notion of an upper rim since it does not express the two 
crucial conditions that $\mathscr Y $ be a sub-class of $\mathscr X$, and that $\mathscr X $ be a positival class and $R$ belong to it. This makes 
`$\mathscr XR \sqcup S^{-}$' do not render the condition that $ S$ be smaller than $ R$ over $\mathscr X$.
Hence, according to (\ref{3.6}$'$), it could happen that  $\left[\left(\mathscr L\right) \Finv \left(\mathscr H\right)\right]\left(\mathbf R\right)$ even if
$\mathscr H $ is not a sub-class of $\mathscr L$ or $\textbf R$ it is not such that any relation
smaller than it over $\mathscr L$ has $\mathscr H$.
It follows that the conjunction
\begin{equation*}
\mathcal{L}\left(\mathscr L \right) \wedge \mathscr L\textbf R \wedge \left[\left(\mathscr L\right) \Finv \, \left(\mathscr H\right)\right] \left(\textbf R\right)
\end{equation*}
renders the informal condition that $\textbf R$ be an upper rim of $\mathscr H$ over $\mathscr L$ except for the 
requirement that $\mathscr H$ be a sub-class of $\mathscr L$.
What are the consequences of missing this requirement? To see it, let us write the implication
 \begin{equation*}
`\forall S \left[\left( \mathscr XS \wedge \mathscr X\textbf R \sqcup S^-\right)\Rightarrow \mathscr YS\right]\text{'} \qquad \mathrm{as} \qquad `\neg\exists S \left[\mathscr XS \wedge \mathscr X\textbf R \sqcup S^- \wedge \neg \mathscr YS \right]\text{'}. 
\end{equation*}
It this clear that this formula can be satisfied by $\mathscr L$  (as value of $\mathscr X$) and $\mathscr H$ (as value of $\mathscr Y$) even if $\mathscr H$ is not a sub-class of $\mathscr L$. For instance, this is just what happens, whatever first-level binary relation $ R$ might be, if $\mathscr L$ is a sub-class of $\mathscr H$. Hence, missing the mentioned requirement results in admitting that, for any $ R$, if  $\mathscr L$ is a sub-class of $\mathscr H$,
then $\left[\left(\mathscr L\right) \Finv \, \left(\mathscr H\right)\right] \left(R\right)$.
But suppose that the first-level binary relations that have $\mathscr L$ form a positival class and that $\mathbf R$ be one of them, that is, that
$\mathcal{L}\left(\mathscr L\right) \wedge \mathscr L \mathbf{R}$.
If $\mathbf R$ is not the smallest relation that has $\mathscr L$, there is certainly another relation 
$S$ such that $\mathscr L  S \wedge \mathscr L  \mathbf R \sqcup  S^-$. Hence, for it to hold that  
\begin{equation*}
\neg\exists S \left[\mathscr L S \wedge \mathscr L \mathbf R \sqcup S^- \wedge \neg \mathscr H S \right] \qquad \mathrm{and} \qquad [\left(\mathscr L\right) \Finv \left( \mathscr H\right)]\left(\mathbf R\right),
\end{equation*}
it is necessary that any such $ S$ have
$\mathscr H$. But if this is so, then $\mathscr L$ and $\mathscr H$ are not disjoint. This having been established, rewrite the right-hand side of 
(\ref{3.7}$'$) in agreement with (\ref{3.6}$'$), i.e. 
as follows
\begin{equation*}
\begin{array}{c} \mathcal{L} \left(\mathscr X\right) \wedge \mathscr X\mathbf R \wedge 
\neg \exists S \left[ \mathscr X S  \wedge \mathscr X \mathbf R \sqcup S^-  \wedge  \neg \mathscr YS \right] \wedge  
\neg \exists T \left[ \mathscr X T  \wedge \mathscr X T \sqcup \mathbf R^-  \wedge 
\neg \exists W \left( \mathscr X W  \wedge \mathscr X T \sqcup W^-  \wedge \neg \mathscr YW \right)\right].
\end{array}
\end{equation*}
For this conjunction to hold, it has to exist a first-order binary relation that has 
$ \mathscr X $ but not  $\mathscr Y$. The case where  $ \mathscr L $ is a sub-class of $ \mathscr H $ is then expunged from those in which it can happen that $[\left(\mathscr L\right) \textrm{\l} \left(\mathscr H\right)]\left(\mathbf R\right)$ for some $\mathbf R$. Insofar as (\ref{3.7}$'$) implies that 
$[\left(\mathscr L \right) \text{\l } \left(\mathscr H\right)]\left(\mathbf{R}\right)$ only if $\mathscr L$ is positival, $\mathbf{R}$ has it, and 
$[\left(\mathscr L \right) \Finv \left( \mathscr H\right)]\left(\mathbf{R}\right)$, it follows that, provided that $\mathbf{R}$ be not the smallest relation having $\mathscr L$, it can happen that  $[\left(\mathscr L \right) \text{\l } \left(\mathscr H\right)]\left(\mathbf{R}\right)$ 
only if $\mathscr H$ is a sub-class of $\mathscr L$ or, at least, $\mathscr L$ and $\mathscr H$ are not disjoint, but  $\mathscr L $ is not a sub-class of 
$ \mathscr H $. Let now $\mathscr P$ be a property of first-level binary relations satisfying the right-hand side of (\ref{3.8}$'$), and, 
then, such that $\mathcal P \left(\mathscr P\right) $. The sub-group formed by the relations having it is, then, dense. If $\mathbf R$ has  $ \mathscr P $, it cannot happen that it be the smallest relation having it. Hence, it can happen that $\left[\left(\mathscr P\right) \Finv \left(\mathscr H\right)\right]\left(\mathbf{R}\right)$, only if  $\mathscr H $ is either a subclass of $\mathscr P $, or  $\mathscr P $ and  $\mathscr H $ are not disjoint, but  $\mathscr P $ is not a subclass of  $\mathscr H $. Hence  $\mathscr P $ and  $\mathscr H $ are not disjoint, that is, some  relation having $\mathscr H $ has also $\mathscr P$. Thus, even if $\mathscr H $ is not a sub-class of $\mathscr P $, it can happen that it has both some upper rims and an upper limit 
over $\mathscr P $, which is just what is relevant for both (\ref{3.7}$'$) and (\ref{3.8}$'$) to comply with the informal explanations given above.}.

From the previous remarks, it should be clear that if the second-level monadic property $\mathscr P$
is such that $\mathcal{P} \left(\mathscr P\right)$, then the permutations that respectively satisfy `$\eth \mathscr P R$' or 
 `$\eth \left( %
\mathscr P\right) \left( R\right) $', if any, form a totally-ordered, dense and
Dedekind-complete group of permutations. This is just what a domain of magnitudes is in our setting.
Indeed, for Frege, domains of magnitudes are nothing but domains of
positive classes. This suggests either the following explicit definition licensed, again, only by an appropriate form of fourth-order comprehension,
\begin{equation}
\forall \mathscr X \left[
\mathcal{M}\mathscr X \Leftrightarrow \exists \mathscr Y\left[ \mathcal{P}%
\mathscr Y \wedge \forall R\left[ \mathscr X R\Leftrightarrow 
\eth \mathscr Y R \right] \right] \right],\label{3.9}
\end{equation}
or the following abbreviation stipulation requiring 
no comprehension stronger than
(PCA$^{2}_{\Delta^{1}_{0}}$),
\begin{equation}
\mathcal{M}\left( \mathscr X\right) :=\exists \mathscr Y\left[ \mathcal{P}%
\left( \mathscr Y\right) \wedge \forall R\left[ \mathscr XR\Leftrightarrow 
\eth \left( \mathscr Y\right) \left( R\right) \right] \right].  \tag{\ref{3.9}$'$}
\end{equation}
where both `$\mathcal{M} \mathscr X$' and `$\mathcal{M}\left( \mathscr X\right)$' are to be read as `$\mathscr X$ is a domain of magnitudes' or, more
precisely, `the first-level binary relations having $\mathscr X$ form a domain of magnitudes'.

\section{Which Definition, in which System?\label{Sect.LE}}
As a matter of fact, (\ref{3.9}) and (\ref{3.9}$'$) provide two different
definitions of domains of magnitudes. The former results from the explicit definition of the third-order predicate constant `$\mathcal{M}$'.
The latter merely exhibits the third-order open formula briefly
designated by `$\mathcal{M}\left( \mathscr X\right) $', and involves no explicit definition other than (Inversion) and (Composition). Both
render Frege's definition, but  require different logical resources, and play distinct roles in our setting. 

Let us begin with the logical resources they require. A first difference is manifest: while (\ref{3.9}) requires a fourth-order system, 
a third-order system is enough for (\ref{3.9}$'$).
Though both systems encompass no proper axioms, the former is, by far, more entangled than the latter. This is not only because of its
higher order, but also because of the forms of comprehension it has to incorporate, in order to license (\ref{3.9}).  
Besides (PCA$^{2}_{\Delta^{1}_{0}}$)---or its instances 
(Composition-CA) and (Inversion-CA)---, required to license (Composition) and (Inversion), it also  
calls for other comprehension axioms, 
respectively required to license (Functionality) and  (\ref{3.4}-\ref{3.9})\footnote{%
Namely:
\begin{equation}
\exists \mathscr X
\forall R \left[\mathscr X R\Leftrightarrow \phi_{\Delta_0^1}\right] \qquad 
\begin{tabular}{l}
[where `$\phi_{\Delta_0^1}$' stands for a second-order predicative formula],
\end{tabular}
\tag{CA$^{3}_{\Delta_0^1}$}
\end{equation}
required to license (Functionality);
\begin{equation}
\forall \mathscr X \exists  \mathscr Y \forall R \left[\mathscr Y R \Leftrightarrow \phi_{\Sigma_{1(1)}^2}\left(\mathscr X\right)\right]
\qquad
\left[
\begin{tabular}{l}
where `$\phi_{\Sigma_{1(1)}^2}\left(\mathscr X\right)$' stands for a third-order formula \\ involving a second-order existential quantifier \\
and the parameter `$\mathscr X$'
\end{tabular}
\right],
\tag{PCA$^{3}_{\Sigma_{1(1)}^2}$}
\end{equation}
required to license (\ref{3.4});
\begin{equation}
\exists \mathcal X \forall \mathscr X \left[\mathcal X \mathscr X \Leftrightarrow \phi_{\Pi_{1(1)}^{2}}\right]
\qquad
\left[
\begin{tabular}{l}
where `$\phi_{\Pi_{1(1)}^{2}}$' stands for a third-order formula \\ involving a
second-order universal quantifier
\end{tabular}
\right],
\tag{CA$^{4}_{\Pi_{1(1)}^{2}}$}
\end{equation}
required to license (\ref{3.2});
\begin{equation}
\exists \mathcal V \forall \mathscr X, \mathscr Y \forall R  \left[\mathscr X \mathcal V_{R} \mathscr Y  
\Leftrightarrow \phi_{\Pi_{1(1)}^{2}}\right] \qquad
\begin{tabular}{l}
[where `$\phi_{\Pi_{1(1)}^{2}}$' is as in (CA$^{4}_{\Pi_{1(1)}^{2}}$)],
\end{tabular}
\tag{CA$^{4}_{\Pi_{1(1)}^{2}}$}
\end{equation}
required to license (\ref{3.6});\newpage
\noindent \begin{equation}
\exists \mathcal V \forall \mathscr X, \mathscr Y \forall R  \left[\mathscr X \mathcal V_{R} \mathscr Y  
\Leftrightarrow \phi_{\Sigma_{1(1)}^{3}}\right]
\qquad
\left[
\begin{tabular}{l}
where  `$\phi_{\Sigma_{1(1)}^{3}}$' stands for a fourth-order formula \\ involving a
second-order existential quantifier
\end{tabular}
\right],
\tag{CA$^{4}_{\Pi_{1(1)}^{3}}$}
\end{equation}
required to license (\ref{3.7});
\begin{equation}
\exists \mathcal X \forall \mathscr X \left[\mathcal X \mathscr X \Leftrightarrow \phi_{\Pi_{1(2)}^{3}}\right]
\qquad
\left[
\begin{tabular}{l}
where  `$\phi_{\Pi_{1(2)}^{3}}$' stands for a fourth-order formula \\ involving a
third-order universal quantifier
\end{tabular}
\right],
\tag{CA$^{4}_{\Pi_{1(2)}^{3}}$}
\end{equation}
required to license (\ref{3.8});
\begin{equation}
\exists \mathcal X \forall \mathscr X \left[\mathcal X \mathscr X \Leftrightarrow \phi_{\Sigma_{1(2)}^{3}}\right]
\qquad
\left[
\begin{tabular}{l}
where  `$\phi_{\Pi_{1(2)}^{3}}$' stands for a fourth-order formula \\ involving a
third-order existential quantifier
\end{tabular}
\right],
\tag{CA$^{4}_{\Pi_{1(2)}^{3}}$}
\end{equation}
required to license (\ref{3.9}).}. The latter system only needs to involve (PCA$^{2}_{\Delta^{1}_{0}}$), or merely (Composition-CA) and (Inversion-CA), instead.

First-order variables (and the usual logical constants) apart, the language of this latter system has only to include 
dyadic second-order and monadic
third-order variables, together with the corresponding quantifiers\footnote{Notice that a third-order quantifier only occurs once: 
in the right-hand side of (\ref{3.8}$'$). This single occurrence is however 
essential for the definition of domains of magnitude to succeed.}, and the functional constants `$^{-}$' and `$\sqcup $'. Though 
third-order, this system is, then, quite weak. As a matter of fact, we have nevertheless shown that Frege's
definition of domains of magnitudes can be consistently rephrased in such a weak system, and is, then, so rephrased,
equiconsistent  with it\footnote{A note of caution. The fact that, when rephrased as suggested, 
Frege's definition is equiconsistent with this system
does not entail at all that the original definition requires no impredicative comprehension, 
and is, then, consistent in itself. It crucially involves the function $\xi \frown \zeta$, which cannot be defined without impredicative (second-order) comprehension.}. For future reference, call this system `L$_2$PCA$^{2}_{\Delta^{1}_{0}}$'.

It remains, however, that (\ref{3.9}$'$), and, then, this very system, are suitable for our present purpose 
only if we are content with
admitting that a second-level property $\mathscr M$
is a domain of magnitudes (or that the first-level binary relations that have
it form such a domain)  if and only if it satisfies the right-hand side of (\ref{3.9}$'$). Were we, instead, in need of defining
a (third-level) property that $\mathscr M$ has if and only if it is so, (\ref{3.9}$'$) would no more be suitable, and 
we would have to recur to (\ref{3.9}) and the corresponding fourth-order and much stronger system. Provided that
the definition of positival and positive classes is, in the present setting, merely instrumental to that of
domains of magnitudes,  the former attitude might be easily admitted for the corresponding properties. But one might consider that
the same attitude is not admissible in the case of these very domains, whose definition is the final outcome of Frege's work, on pain of missing a genuine entity counting as such a domain, and, then, the definition itself.

Still, even if the definition of domains of magnitudes is the last step Frege reached in his formalization
of real analysis\footnote{To be more precise, Frege offers
no explicit formal definition of domains of magnitudes, and rather is
content with  informally claiming that a domain of magnitudes is the domain of
a positive class. \label{ft.TBMP}}, it is in no way the final step such a formalization should have reached: this should have rather
been an explicit definition of real numbers, and, possibly, of the operations and relations making them form a 
totally ordered and Dedekind-complete field.
Hence, if Frege's informal indications for reaching this final aim may be rendered in
our setting without defining the predicate constant `$\mathcal{M} $', there is no
stringent reason for accepting the foregoing argument, so that (\ref{3.9}$'$) and L$_2$PCA$^{2}_{\Delta^{1}_{0}}$
may be  considered enough for the purpose of rendering the result he was envisaging.
In \S~\ref{Sect.RBOEP}, we will show that this can be actually done. We can, then, conclude that, whereas
(\ref{3.9}), and the fourth-order system it requires, are suitable for the purpose of defining domains of magnitude, as such,
(\ref{3.9}$'$) and  L$_2$PCA$^{2}_{\Delta^{1}_{0}}$ are so for the purpose of defining real numbers as ratios on such domains. 
As this is our goal, we'll go, then, for the latter option. 

This is all the more justified because no form of comprehension can
guarantee the
existence of positival and positive classes and domains thereof. Surely, by the standard interpretation of higher-order logic, the stipulations (\ref{3.2}$'$) and (\ref{3.8}$'$-\ref{3.9}$'$) being given, the following
instances of third-order impredicative comprehension
\begin{equation*}
\begin{array}{l}
\exists \mathscr X \forall R \left[ \mathscr X R \Leftrightarrow \exists %
\mathscr Y \left( \mathcal{L} \left( \mathscr Y \right) \wedge \mathscr Y R
\right)\right], \medskip \\
\exists \mathscr X \forall R \left[ \mathscr X R \Leftrightarrow \exists %
\mathscr Y \left( \mathcal{P }\left( \mathscr Y \right) \wedge \mathscr Y R
\right)\right], \medskip \\
\exists \mathscr X \forall R \left[ \mathscr X R \Leftrightarrow \exists %
\mathscr Y \left( \mathcal{M }\left( \mathscr Y \right) \wedge \mathscr Y R
\right)\right]
\end{array}
\end{equation*}%
are enough to ensure the existence of the second-level properties that a
first-level binary relation has to have for being respectively included in a
positival and a positive class and in a domain of magnitudes.
Again, the following instances of fourth-order predicative comprehension 
\begin{equation*}
\begin{array}{l}
\exists \mathcal{X }\forall \mathscr X \left[ \mathcal{X }\mathscr X
\Leftrightarrow \mathcal{L} \left( \mathscr X\right) \right],
\medskip \\
\exists \mathcal{X }\forall \mathscr X \left[ \mathcal{X }\mathscr X
\Leftrightarrow \mathcal{P }\left( \mathscr X\right)\right],
\medskip \\
\exists \mathcal{X }\forall \mathscr X \left[ \mathcal{X }\mathscr X
\Leftrightarrow \mathcal{M }\left( \mathscr X\right)\right]
\end{array}
\end{equation*}%
are enough to ensure the existence of the third-level properties
of being a positival and a positive class and a domain of magnitudes.
Still, securing this existence is no substantial
achievement, since these properties would exist even if they were empty, or there were
not enough relations satisfying them, to make them play the required role 
in a definition of the reals.

Even more so, if we grant the extensional conception of
properties and relations, we also have to grant that, no matter how the first-order
domain might be, both the empty first-level binary relation and the empty second-level
property exist\footnote{Frege is not clear about what 
he takes to make a property or a relation exist, if at all. It is even plausible to ascribe
him an intensional conception, which makes the talk of existence of properties and
relations inappropriate (Panza 2015). What is unquestionable is that he does not explicitly endorse any sort of comprehension axiom, 
by rather pervasively admitting of a substitution policy which we could consider as 
equivalent to full second-order 
comprehension. In the light of BLV, our replacing value-ranges of 
binary first-level relations with these very relations seems, however, in line with granting the extensional conception 
of properties and relations while doing semantic considerations
about L$_2$PCA$^{2}_{\Delta^{1}_{0}}$. 
This means informally taking a $n+1$-order $m$-adic predicate variable to range on all the subsets of ordered $m$-tuples 
of the elements of the  $n$-order domain  ($n,m=1,2,\ldots$), and a predicate constant to designate one of these subsets.}, 
which alone is enough to 
ensure that all six foregoing
second- and third-level properties
exist, in turn, and are nonempty, since, according to (\ref{3.2}$'$) and (\ref{3.8}$'$-\ref{3.9}$'$),  the empty second-level property is
at the same time a positival class, a positive one, and a domain of
magnitudes\footnote{Let $\mathscr E_{0}$ be the second-level empty property. From (\ref{3.2}$'$) it follows that $\mathcal L (\mathscr E_{0})$, since $\forall R\lnot\left[\mathscr E_{0}R\right]$. Let $\mathsf{V}$ 
be the empty first-level binary relation, and 
$\mathscr E_{1}$ the second-level property of being this property. If $\mathbf R$ is a first-level binary relation (extensionally) 
distinct from $\mathsf{V}$, then, we have that both $\mathscr E_{1} \mathsf{V}$ and  $\lnot \mathscr E_{1} \mathbf R$.
Moreover, from (Funcionality), (Inversion$'$) and (Composition$'$), it follows that $\mathscr I(\mathsf{V})$, and that $\mathsf{V}$
extensionally  coincides both with $ \mathsf{V}^{-}$ and with $ \mathsf{V} \sqcup  \mathsf{V}^{-}$, and, then, also with
$ \mathsf{V}^{-} \sqcup  \mathsf{V}$ (see Frege, 1893-1903, \S~II.164).
This is enough to show that $\mathsf V$ does belong to no positival (and, then, positive) class, and that 
both $\lnot \mathcal L (\mathscr E_{1})$ and 
$ \mathcal L^{*} (\mathscr E_{1})$, where `$\mathcal L^{*}(\mathscr X)$' is the abbreviation of the formula resulting
from the right-hand side of (\ref{3.2}$'$), by skipping the conjunct `$\lnot \mathscr X   R  \sqcup  R^{-}$'.
Then, $\mathscr E_{1}$ is positival, except for admitting the possibility that the identity relation be included in it. 
Thus, both the second-level properties  $\left[ R : \exists \mathscr Y \left( \mathcal L \left( \mathscr Y\right) \wedge \mathscr Y R \right)\right]$ and  $\left[ R : \exists \mathscr Y \left( \mathcal L^{*} \left( \mathscr Y\right) \wedge \mathscr Y R \right)\right]$ and the third level
ones $\left[ \mathscr X : \mathcal L \left( \mathscr X \right)\right]$ and $\left[ \mathscr X : \mathcal L ^{*}\left( \mathscr X \right)\right]$
not only exist by appropriate forms of comprehension, but are also nonempty. 
Consider now (\ref{3.8}$'$). It is enough to observe that $\mathcal L (\mathscr E_{0})$ and $\lnot \exists R \left[ \mathscr E_{0} R\right]$ to conclude that $\mathcal P (\mathscr E_{0})$. Look, then, at $\mathscr E_{1}$. 
Insofar as $\mathbf R \sqcup  \mathsf{V}^{-}$ coincides with $\mathsf V$, whatever first-level binary relation $\mathbf R$ might be, from 
(\ref{3.6}$'$) it follows that $\left[ (\mathscr E_{1}) \Finv  (\mathscr Y) \right](\mathbf R)$ if and only if $\mathscr Y$ is $\mathscr E_{1}$ itself. Insofar as $\forall S \lnot \left[ \mathscr E_{1} S \wedge \lnot \mathscr E_{1} S\right]$,
from (\ref{3.8}$'$) it also follows that $\mathcal P^{*}(\mathscr E_{1})$, where `$\mathcal P^{*}(\mathscr X)$' is the abbreviation of the formula resulting
from the right-hand member of (\ref{3.8}$'$) by replacing
`$\mathcal L \left(\mathscr X\right)$' with `$\mathcal L^{*} \left(\mathscr X\right)$'. So, the properties $\left[ R : \exists \mathscr Y \left( \mathcal P \left( \mathscr Y\right) \wedge \mathscr Y R \right)\right]$,  $\left[ R : \exists \mathscr Y \left( \mathcal P^{*} \left( \mathscr Y\right) \wedge \mathscr Y R \right)\right]$, $\left[ \mathscr X : \mathcal P \left( \mathscr X \right)\right]$ and $\left[ \mathscr X : \mathcal P ^{*}\left( \mathscr X \right)\right]$ not only exist by appropriate forms of comprehension, but they are also nonempty. Finally,
consider (\ref{3.9}$'$).
Provided that `$\mathcal M^{*}(\mathscr X)$' be the abbreviation of the formula resulting
from the right-hand side of (\ref{3.9}$'$) by replacing
`$\mathcal P \left(\mathscr X\right)$' with `$\mathcal P^{*} \left(\mathscr X\right)$', it is immediate that both $\mathcal M (\mathscr E_{0})$ and $\mathcal M^{*} (\mathscr E_{1})$,  just because $\mathcal P (\mathscr E_{0})$ and $\mathcal P^{*} (\mathscr E_{1})$ and the domains of  $\mathscr E_{0}$ and $\mathscr E_{1}$ respectively coincide with $\mathscr E_{0}$ and $\mathscr E_{1}$ themselves. 
Hence, the properties $\left[ R : \exists \mathscr Y \left( \mathcal P \left( \mathscr Y\right) \wedge \mathscr Y R \right)\right]$,  $\left[ R : \exists \mathscr Y \left( \mathcal P^{*} \left( \mathscr Y\right) \wedge \mathscr Y R \right)\right]$, $\left[ \mathscr X : \mathcal M \left( \mathscr X \right)\right]$ and $\left[ \mathscr X : \mathcal M ^{*}\left( \mathscr X \right)\right]$, too, not only exist by appropriate forms of comprehension, but are also nonempty. Notice that $\mathcal M ^{*}$ does not (extensionally) coincide with 
$\mathcal M$, since
$\mathcal M^{*} (\mathscr E_{1})$, but not $\mathcal M (\mathscr E_{1})$.
\label{CVDM}}.

Clearly though, if a positival class exists that
includes at least a first-level binary relation, then it necessarily
includes infinitely many (extensionally) distinct relations, which are continuously many if
the class is positive\footnote{
More precisely, this is with respect to the full model of  L$_2$PCA$^{2}_{\Delta^{1}_{0}}$,
which we also take as its intended one,
where the first-order variables, $x$, $y$, $z$, \ldots vary over a large enough domain $\mathbf A$ of objects,
the second-order dyadic variables $R$, $S$, $T$, \ldots vary over the full power set of  $\mathbf A \times \mathbf A$, and 
the third-order
monadic variables $\mathscr X$, $\mathscr Y$, \ldots vary over the full power set of the power set of  $\mathbf A$.
It is, indeed, easy to see that L$_2$PCA$^{2}_{\Delta^{1}_{0}}$ is far from categorical. A simple way to see it (thanks to Andrew Moshier for suggesting it to us) is to observe that the lack
of third-order comprehension makes this system have a model where
the third-order variables vary over the empty set. In this model, all closed formulas
beginning by a third-order universal quantifier, as the third conjunct in the right-hand side of (\ref{3.8}$'$), are vacuously true. 
This makes `$\mathcal P\left(\mathscr X \right)$' trivially satisfied by countably many (appropriate) binary 
relations (whereas PCA$_{\Delta_{0}^{1}}^{2}$ makes any model of L$_2$PCA$^{2}_{\Delta^{1}_{0}}$
include countably many first-level binary relations). 
One might be surprised we take the full model to be the intended model, rather than an appropriate Henkin one.
Still, we think a restriction on comprehension, or even the lack of it, is no reason 
for imposing a restricted semantics: one thing is
the logical resources, in particular the instances of comprehension, required for a definition; another 
the selection of the intended model. The former are deductive, syntactical tools required to formulate definitions; the 
latter depends on semantic considerations relative to the informal piece of knowledge that the relevant system is
expected to render---which, in this case, involves the idea of a positive class as a complete semi-group of permutations.
We are not going to take a stand on this matter, here, but costs and benefits are worth mentioning. Should our definition be required to recover positive classes including 
exactly 2$^{\aleph_0}$ permutations, then the intended semantic is to be clearly full.
If impredicativity were deemed a price worth paying for this purpose, one should also be aware of the heavy-duty logical resources called for it. \label{ft.catMod}}.
And, if this holds for a positive class, it also
holds for its domain. It follows that a positive class and a domain of
magnitudes exist, which include at least a first-level binary relation, if
and only if there exists one including continuously many of them.
No semantic consideration is, however, apt to prove the existence of nonempty
positive classes and domains thereof. To prove it, it is
necessary to prove that there are enough objects, for defining on them
continuously many appropriate relations.

This was the crucial challenge for Frege's definition of the real numbers,
and it is also ours. Before tackling it, it is, however, in order to see what makes this proof indispensable
both for Frege's and for our purposes. This requires looking at
how Frege's original theory of the reals as objects can be revived in our setting, on the basis of
(\ref{3.8}$'$).

\section{Real Numbers as Objects\label{Sect.RBOEP}}
Apart from some generalities, most of which we
already discussed above, and the sketchy outline of a possible
existence proof for nonempty domains of magnitudes, which we will consider
in \S ~\ref{Sect.BP}, the \textit{pars construens} of part III of 
\textit{Grundgesetze} only contains the definition of these domains. No 
precise indication of how to define the real numbers is available. The only thing that is clear is that these numbers should
be defined as ratios of magnitudes, and that these ratios have to be defined as objects, i.e.,
first-order items, logically speaking\footnote{From what he writes in the very \S\ of his treatise, it seems that Frege was also
requiring that real numbers  form themselves a domain of magnitudes (Frege 1893-1903, \S~II.245;
Frege 2013, p. 243$_2$):
\begin{quote}
The commutative law for the
domain of a positive class is thus proven.
The next task is now to show that there
is a positive class, as indicated in \S~164 [see p.~\pageref{164}, above].
This opens the possibility of defining a
real number as a ratio of magnitudes
of a domain that belongs to a positive
class. Moreover, we will then be able to
prove that the real numbers themselves
belong as magnitudes to the domain of
a positive class.
\end{quote}
This further requirement would have not only uselessly entangled Frege's own first-order definition, if he 
completed it (Dummett 1991, pp.~190-91), but it is also logically incompatible, in 
our predicate setting, with the requirement that real numbers be objects.
This is why we will set it aside in  our reconstruction.}.
 
\subsection{Euclid's Principle with Natural Numbers}
A simple way to accomplish this plan is by
an abstraction
principle governing an operator taking pairs of relations (i.e. magnitudes) from a
domain of magnitudes as arguments and having objects (i.e. ratios) as values. As
suggested by Simons (1987, pp. 39--40) and Dummett (1991, pp. 290--91),
this can be done by rephrasing definition V.5 of Euclid's \textit{Elements}, and defining the relation of
proportionality between four magnitudes taken two by two\footnote{%
As observed by Dummett (1991, pp. 282--83), Euclid's definition, probably tracing back to Eudoxus, in fact,
had been explicitly appealed to by H\"{o}lder in his paper on the
\textquotedblleft axioms of quantity\textquotedblright\ (H\"{o}lder 1901),
which appeared two years before the second volume of the \textit{Grundgesetze%
}. But, apparently, Frege's was not aware of this. On H\"{o}lder (1901), cf. B\l{}aszczyk (2013), which also sums up how the notion of magnitude was investigated around the end of 19th century by several mathematicians, including Du Bois-Reymond, Stolz, and Weber, by explicitly referring to Euclid's theory, and achieving results mathematically equivalent to Frege's.}

This raises a technical difficulty, however. Whereas this definition applies
only if the magnitudes composing each pair are such that it cannot be the case
that any multiple of
one of them be smaller, equal, or greater than any multiple of the other, 
this condition is not
met by pairs of magnitudes of the same domain, since, differently
from what happened for Euclid's ones, these are intended to be either
positive, negative, or null.
A way to solve the difficulty is by appropriately modifying the very structure of Euclid's definition, in order to get the following
abstraction principle, which deserves, nevertheless, the name of `Euclid's Principle' (or `EP', from now on),
\begin{equation*}
\left\{ 
\begin{array}{l}
\forall _{\left( \mathcal{P}\right) }\mathscr X,\mathscr X^{\prime}\medskip \\ 
\forall _{\left( \eth \left( \mathscr X\right) \right) }R\,\forall _{\mathscr %
X}S\medskip \\ 
\forall _{\left( \eth \left( \mathscr X^{\prime}\right) \right) }R^{\prime }\,\forall _{%
\mathscr X^{\prime}}S^{\prime }%
\end{array}%
\right\} \left[ 
\begin{array}{l}
\mathfrak{R}\left[ R,S\right] =\mathfrak{R}\left[ R^{\prime },S^{\prime }%
\right] \Leftrightarrow \medskip \\ 
~~\left[ 
\begin{array}{l}
\left[ 
\begin{array}{l}
\mathscr XR\wedge \mathscr X^{\prime}R^{\prime }\wedge \medskip \\ 
\forall _{\mathbf{N}}x,y(xR\sqsubset_{\mathscr X}yS\Leftrightarrow xR^{\prime }\sqsubset_{%
\mathscr X^{\prime}}yS^{\prime })%
\end{array}%
\right] \vee \medskip \\ 
\left[ 
\begin{array}{l}
\mathscr XR^{-}\wedge \mathscr X^{\prime}R^{\prime -}\wedge \medskip \\ 
\forall _{\mathbf{N}}x,y(xR^{-}\sqsubset_{\mathscr X}yS\Leftrightarrow xR^{\prime
-}\sqsubset_{\mathscr X^{\prime}}yS^{\prime })%
\end{array}%
\right] \vee \medskip \\ 
\left[ \forall z,w\left[ 
\begin{array}{l}
\left( zRw\Leftrightarrow z[S\sqcup S^{-}]w\right) \wedge \medskip \\ 
\left( zR^{\prime }w\Leftrightarrow z[S^{\prime }\sqcup S^{\prime -}]w\right)%
\end{array}%
\right] \right]%
\end{array}%
\right]%
\end{array}%
\right], \tag{EP} 
\end{equation*}%
\noindent where:
\begin{itemize}
 \item[--] `$\mathfrak{R}$' is the relevant abstraction operator;
\item[--] `$\forall _{\left( \mathcal{P}\right) }\mathscr X\left[ \varphi \right]
$' abbreviates `$\forall \mathscr X\left[ \mathcal{P}\left( \mathscr %
X\right) \Rightarrow \varphi \right] $', and the same for `$\mathscr X^{\prime}$';
\item[--] `$\forall _{\left( \eth \left( \mathscr X\right) \right) }R\left[
\varphi \right] $' abbreviates `$\forall R\left[ \eth \left( \mathscr %
X\right) \left( R\right) \Rightarrow \varphi \right] $', and the same for `$\mathscr X^{\prime}$' and `$%
R^{\prime }$';
\item[--] `$\forall _{\mathscr X}S\left[ \varphi \right] $' abbreviates `$%
\forall S\left[ \mathscr XR\Rightarrow \varphi \right] $', and the same for `$\mathscr X^{\prime}$' and `$S^{\prime }$';
 \item[--] `$\forall _{\mathbf{N}}x\left[ \varphi \right] $' abbreviates `$%
\forall x\left[ \mathbf{N}x\Rightarrow \varphi \right] $', and the same for `%
$y$';
\item[--] `$\mathbf{N}$' denotes the property of being a natural number;
 \item[--] `$xR$' abbreviates `$\underbrace{R\sqcup R\sqcup \ldots \sqcup R}%
_{x~{\text{times}}}$, and the same for `$y$', `$S$',  `$S^{\prime}$', $R^{\prime }$', `$R^{-}$', `$%
R^{\prime -}$';
 \item[--] `$xR\sqsubset_{\mathscr X}yS$' abbreviates `$\mathscr XyS\sqcup (xR)^{-}$', and the same for `$\mathscr X^{\prime}$', `$R^{\prime }$', `$S^{\prime }$', `$R^{-}$%
', `$R^{\prime -}$'.
\end{itemize}

Informally speaking, EP states that for whatever pairs of domains of 
magnitudes, issued by two positive classes $\mathscr P$ and $\mathscr P^{\prime}$\footnote{Notice that EP 
does not involve domains of magnitudes as such, but rather positive classes and their domains.
This is perfectly in line with Frege's missing a formal definition of these domains: see footnote (\ref{ft.TBMP}), above.}, and whatever
two ordered pairs of permutations $\mathbf R$, $\mathbf S$ and $\mathbf R'$, $\mathbf S'$, such that
$\mathbf R$, and $\mathbf R'$ respectively belong to the domains of these classes, while 
$\mathbf S$ and $\mathbf S'$ belong to the classes themselves (and are, then, intended as positive), the ratio  $\mathfrak{R}\left[ \mathbf R,\mathbf S\right]$
of the elements of the
first pair is the same as the ratio  $\mathfrak{R}\left[ \mathbf R',\mathbf S'\right]$
of the elements of the second pair if
and only if:

\begin{itemize}
\item[--] either both the first elements $\mathbf R$ and $\mathbf R'$ of these pairs belong to the respective
positive classes $\mathscr P$ and $\mathscr P^{\prime}$ (and are, then, intended as positive), and their equimultiples are always  smaller than the equimultiples of the second elements\footnote{Notice also
that, whereas Euclid's definition requires that the equimultiples of the first and the third, among the four
relevant magnitudes, ``alike exceed, are alike equal to, or alike fall short of 
[{\greektext{<'ama <uper'eqh| >`h <'ama >'isa  \~>h| >`h <'ama >elle'iph|}}]'' (Euclid EH, vol.~II, p.~114)
the  equimultiples of the second and the fourth, we can just 
require that the equimultiples of the first and the third magnitudes all be smaller than those of the second and the fourth, since, as noticed by Scott (1958-59),
in the case of Archimedean magnitudes, the latter condition entails the former.};

\item[--] or both the inverses $\mathbf R^{-}$ and $\mathbf R'^{-}$ of the first elements $\mathbf R$ and $\mathbf R'$ of these pairs belong to the respective positive classes $\mathscr P$ and $\mathscr P^{\prime}$ (so that 
$\mathbf R$ and $\mathbf R'$ are intended as negative), and their equimultiples are always
smaller than the equimultiples of the second elements;

\item[--] or both the first elements $\mathbf R$ and $\mathbf R'$  of these pairs are the identity relation (and are, then, intended as null).
\end{itemize}

So rephrasing Euclid's definition surely solves the technical difficulty, but
it does not solve all problems: though EP involves 
neither a predicate constant `$\mathcal P $' for the third-level property of being a positive class,
nor a functional constant `$\eth $' for the domain operator, but merely depends on the stipulations
(\ref{3.8}$'$) and  (\ref{3.4}$'$),
it cannot, as such, be added to L$_2$PCA$^{2}_{\Delta^{1}_{0}}$ as a new axiom, so as to get an extended system in which
real numbers are to be defined. There are two reasons for that. First of all,
EP involves the predicate constant `$\mathbf N$' for the first-level property of
being a natural number, which is not and
cannot be defined within L$_2$PCA$^{2}_{\Delta^{1}_{0}}$.
Secondly, it involves the symbol `$xR$' (or `$\underbrace{R\sqcup R \sqcup \ldots \sqcup R}%
_{x~\text{times}}$') (where `$x$' is a variable ranging over the natural numbers) whose use in a formal system
 is licensed only if this latter contains a device to count the iterated applications of the functional constant `$\sqcup $', which is not 
 and cannot
be provided within L$_2$PCA$^{2}_{\Delta^{1}_{0}}$.

A way to overcome the first issue is to extend
L$_2$PCA$^{2}_{\Delta^{1}_{0}}$  to a stronger system, in which the property of being a natural number can somehow be defined, 
e.g. by adding HP as a new axiom and appropriately extending its language by monadic first-order predicates,
to make it include
Frege Arithmetic (FA). This would be, however, a quite radical move, which would also openly conflict
with Frege's requirement of non-arithmeticity for his definition of real numbers. Even more so, it would not
overcome the second issue, unless the new system were supplied by some ingenious device not usually available in (the current
versions of) FA.

A much less radical and costly, though a bit
laborious, way to overcome both issues at once is available. It is in fact suggested by a 
trick Frege appeals to in proving the
Archimedeanicity of positive classes (Frege, 1893-1903, \S\S~II.199-214). It consists in amending EP with the help of some new abbreviation stipulations, which merely require
adding new third-order binary variables.

\subsection{Euclid's Principle without Natural Numbers\label{Sect.EPWNN}}

Let us begin by adopting the following
new abbreviation stipulation:
\sloppy \[
\mathscr D_{\left(T\right)}\left(R, S\right) := \forall x,y \left(xSy
\Leftrightarrow x[T\sqcup R]y\right).
\]
Let $\mathbf R$, $\textbf S$,
and $\mathbf T$, be whatever first-level binary relations. According to this stipulation,
the formula `$\mathscr D_{\left( \mathbf T\right) }\left( \mathbf R,\textbf S\right) $' asserts
that $\textbf S$ results from, or extensionally coincide with, the composition of $\mathbf T$ and $\mathbf R$. Hence,
 `$\mathscr D_{\left( \mathbf R\right) }\left( \mathbf R,\textbf S\right) $'  asserts that
 $\textbf S$ results from the composition of $\mathbf R$ with itself. In the notation employed in stating EP, this means that $\textbf S$ coincides with $2\mathbf R$. 

This allows to simulate the usual definition of the weak ancestral of a binary relation:
\[
 \mathscr D_{\left( R\right)}^{\widen}\left(S\right) 
:= \forall \mathscr %
X \left[\left[%
\begin{array}{l}
\mathscr X  R \, \wedge \medskip \\ 
\forall R^{\prime }, S^{\prime }\left[ \left[%
\begin{array}{l}
\mathscr X R^{\prime }\wedge \medskip \\ 
\mathscr D_{\left( R^{\prime }\right)}\left( R^{\prime }, S^{\prime }\right)%
\end{array}
\right] \Rightarrow \mathscr X  S^{\prime }\right]%
\end{array}
\right] \Rightarrow \mathscr X  S \right].
\]
This makes the formula `$\mathscr D_{\left( \mathbf R\right)}^{\widen
}\left(\mathbf S\right)$' assert 
that $\textbf S$ extensionally coincides with $\textbf R$ or with the relation resulting from a iterated composition of $\mathbf R$ with 
itself, and is, then, a multiple of $\mathbf R$ itself. In the notation employed in stating EP, this means that $\textbf S$ coincides with $n\mathbf R$, for some natural number  $n$.
Let, now, $\mathscr P$ be a positive class, and $\textbf R$ a relation in it. It is clear that if
$\mathscr D_{\left( \mathbf R\right)}^{\widen
}\left(\mathbf T\right)$ and $\mathscr D_{\left( \mathbf R\right)}^{\widen
}\left(\mathbf S\right)$, then both $\mathbf T$ and $\mathbf S$  belong to  $\mathscr P$. Suppose it is so, and that
 $\mathscr P \mathbf T \sqcup \mathbf S^{-}$. We can, then, take $\mathbf S$ to be smaller than $\mathbf T$ over $\mathscr P$. Hence, if also $\mathscr P^{\prime}$ is a positive
class (either distinct from $\mathscr P$ or not), $\textbf R^{\prime}$ is a first-level binary relation that belongs to it, and it is also the
case that $\mathscr D_{\left( \mathbf R^{\prime}\right)}^{\widen
}\left(\mathbf T^{\prime}\right)$ for some first-level binary relation $\mathbf T^{\prime}$, then  
$\textbf T$ is the same multiple of $\mathbf R$ over
$\mathscr P$ as $\textbf T^{\prime }$ of $\mathbf R^{\prime }$ over $\mathscr P^{\prime }$ if and only if there are
as many first-level binary relations that satisfy the open formula
`$\mathscr D_{\left( \mathbf R\right)}^{\widen
}\left( S\right) \wedge \mathscr P \mathbf T \sqcup S^{-}$' as those that satisfy the other open formula
`$\mathscr D_{\left( \mathbf R^{\prime}\right)}^{\widen
}\left( S^{\prime}\right) \wedge \mathscr P^{\prime} \mathbf T^{\prime} \sqcup S^{\prime-}$'.

This suggests enriching the language of L$_2$PCA$^{2}_{\Delta^{1}_{0}}$ by introducing third-order binary
variables, ranging over second-level binary relations between first-level such relations, and adopting the following further abbreviation stipulation
\begin{equation}
_{\left( \mathscr X,\mathscr X^{\prime }\right) }\mathcal{E}\left(
 R, T, R^{\prime }, T^{\prime }\right) :=
 \left\{
 \begin{array}{l}
 \mathscr D_{\left( R\right)}^{\widen}\left(T\right) \wedge
 \mathscr D_{\left( R^{\prime}\right)}^{\widen}\left(T^{\prime}\right)
 \wedge \medskip \\
 \exists \mathscr R
 \left[
\begin{array}{l}
\forall S \left[ 
\begin{array}{l}%
\left(
\mathscr D_{\left( R\right)}^{\widen}\left(S\right)
\wedge 
S \sqsubseteq _{\mathscr X} T
\right)
\Rightarrow
\medskip \\
\exists ! S^{\prime } 
\left[
S \mathscr R S^{\prime} \wedge 
\mathscr D_{\left( R\right)}^{\widen}\left(S^{\prime}\right)
\wedge 
S^{\prime} \sqsubseteq _{\mathscr X^{\prime}} T^{\prime}
\right]
\end{array}%
\right] \wedge \medskip \\ 
\forall S^{\prime} \left[ 
\begin{array}{l}
\left(
\mathscr D_{\left( R^{\prime}\right)}^{\widen}\left(S^{\prime}\right)
\wedge 
S^{\prime} \sqsubseteq _{\mathscr X^{\prime}} T^{\prime}
\right)
\Rightarrow
\medskip \\
\exists ! S 
\left[
S \mathscr R S^{\prime} \wedge 
\mathscr D_{\left( R\right)}^{\widen}\left(S\right)
\wedge 
S \sqsubseteq _{\mathscr X} T
\right]
\end{array}%
\right]
\end{array}
\right]
\end{array}
\right\}, \label{eq.A}
\end{equation}
where `$\mathscr R $' is such a variable, and
`$S \sqsubseteq _{\mathscr X} T$' abbreviates `$\mathscr X T \sqcup  S^{-} \vee
\forall x,y\left[xSy \Leftrightarrow xTy \right]$', and the same as for `$\mathscr X^{\prime}$',
`$T^{\prime}$'  and `$S^{\prime}$'. Thus, if $\mathscr P$, $\mathscr P^{\prime}$, $\textbf R$, $\textbf R^{\prime}$,  $\textbf T$, and $\textbf T^{\prime}$ 
are as above, then the formula `$_{\left( \mathscr P,\mathscr P^{\prime }\right) }\mathcal{E}\left(
 \mathbf R, \mathbf T, \mathbf R^{\prime }, \mathbf  T^{\prime }\right) $' asserts that 
$\textbf T$ is the same multiple of $\mathbf R$ over 
$\mathscr P$ as $\textbf T^{\prime }$ of $\mathbf R^{\prime }$ over $\mathscr P^{\prime }$.

For short, let us, now, adopt this other abbreviation stipulation:
\[
_{\left(\mathscr X, \mathscr X^{\prime}\right)}\mathcal E^{(R,T,R^{\prime},T^{\prime})}_{(S,U, S^{\prime},U^{\prime})}:=\,
_{\left( \mathscr X,\mathscr X^{\prime}\right) }\mathcal{E}\left(R,T,R^{\prime
},T^{\prime }\right) \wedge \, 
_{\left( \mathscr X,\mathscr X^{\prime}\right) }\mathcal{E}\left( S,U,S^{\prime
},U^{\prime }\right).
\]
EP can, then, be restated as follows:
\begin{equation}
\left\{ 
\begin{array}{l}
\forall _{\left( \mathcal{P}\right) }\mathscr X,\mathscr X^{\prime}\medskip \\ 
\forall _{\left( \eth \left( \mathscr X\right) \right) }R \forall _{%
\mathscr X}S \medskip \\ 
\forall _{\left( \eth \left( \mathscr X^{\prime}\right) \right) }R^{\prime
}\forall _{\mathscr X^{\prime}}S^{\prime }%
\end{array}%
\right\} \left[ 
\begin{array}{l}
\mathfrak{R}\left[ R,S\right] =\mathfrak{R}\left[ R^{\prime },S^{\prime }%
\right] \Leftrightarrow \medskip \\ 
\begin{array}{l}
\left[ 
\begin{array}{l}
\left[ 
\begin{array}{l}
\mathscr X R \wedge \mathscr X^{\prime} R^{\prime} \, \wedge \medskip \\
\forall T, U, T^{\prime}, U^{\prime}
\left[\begin{array}{l}
_{\left(\mathscr X, \mathscr X^{\prime}\right)}\mathcal E^{(R,T,R^{\prime},T^{\prime})}_{(S,U, S^{\prime},U^{\prime})}
\Rightarrow 
\medskip \\
\left(T\sqsubset_{\mathscr X}U\Leftrightarrow T^{\prime }\sqsubset_{\mathscr %
X^{\prime}}U^{\prime }\right)
\end{array}
\right]
\end{array}%
\right] \vee \medskip \\ 
\left[ 
\begin{array}{l}
\mathscr X R^{-} \wedge \mathscr X^{\prime} R^{\prime-} \wedge \medskip \\
\forall T, U, T^{\prime}, U^{\prime}
\left[\begin{array}{l}
_{\left(\mathscr X, \mathscr X^{\prime}\right)}\mathcal E^{(R^{-},T,R^{\prime-},T^{\prime})}_{(S,U, S^{\prime},U^{\prime})}
\Rightarrow 
\medskip \\
\left(T\sqsubset_{\mathscr X}U\Leftrightarrow T^{\prime }\sqsubset_{\mathscr %
X^{\prime}}U^{\prime }\right)
\end{array}
\right]
\end{array}%
\right] \vee \medskip \\ 
\left[ 
\begin{array}{l}
\forall zw\left[ zRw\Leftrightarrow z[S\sqcup S^{-}]w\right] \wedge \medskip
\\ 
\forall zw\left[ zR^{\prime }w\Leftrightarrow z[S^{\prime }\sqcup S^{\prime
-}]w\right]%
\end{array}%
\right]%
\end{array}%
\right]%
\end{array}%
\end{array}%
\right].\tag{EP$^*$}
\end{equation}

It should be easy to verify that, informally speaking, EP$^*$ has the same 
content as EP. But it expresses this content in the language of 
L$_2$PCA$^{2}_{\Delta^{1}_{0}}$, merely enriched by the addition of
third-order binary variables as `$\mathscr R$'. This addition being admitted, EP$^*$ can, then, be added to
this system as a supplementary axiom. Since 
EP$^*$ is an abstraction principle, its left-hand side is a first-order identity (i.e. `$\mathfrak{R}\left[ R,S\right]$' and 
`$\mathfrak{R}\left[ R^{\prime },S^{\prime }\right]$'  are singular terms). Moreover, its right-hand side
involves no constant other than `$^{-}$' and `$\sqcup $'. Hence,
adding it to 
L$_2$PCA$^{2}_{\Delta^{1}_{0}}$ requires no
further comprehension axiom. The theory obtained is, then,
a third-order one, including first-order, second-order binary, and third-order
monadic and binary variables, but only admitting predicative second-order comprehension.

\subsection{Real Numbers\label{RN}}
Though EP$^*$ supplies the required grounds for
defining real numbers as objects, this theory falls short of achieving it. All that
one can do, in the light of it, is informally (and meta-theoretically) identify these
numbers with ratios like  $\mathfrak{R}\left[ \mathbf R,\mathbf  S\right] $. If 
a predicate constant designating the first-level property of being a real number is to be available,
one also has to admit a new form of comprehension, for licensing the following explicit definition:
\begin{equation}
\forall x\left[ \mathsf{R}x\Leftrightarrow \exists _{\left( \mathcal{P}%
\right) }\mathscr X\exists _{\left( \eth \left( \mathscr X\right) \right)
}R\exists _{\mathscr X}S\left[ x=\mathfrak{R}\left[ R,S\right] \right] %
\right].\label{5.1}
\end{equation}
What we need, then, is the following second-order third-orderly impredicative axiom-scheme:
\begin{equation}
\exists X \forall x\left[ Xx\Leftrightarrow \phi_{\Sigma^{2}_{1}} \right] \tag{CA$^{2}_{\Sigma^{2}_{1}}$},
\end{equation}%
where `$\phi_{\Sigma^{2}_{1}}$' stands for a third-order formula involving a third-order  existential
quantifier---together with a second-order one.

It is then only in such a (highly) impredicative third-order theory obtained from L$_2$PCA$^{2}_{\Delta^{1}_{0}}$ by adding to it both the proper axiom EP$^*$ and the
comprehension axiom-schema (CA$^{2}_{\Sigma^{2}_{1}}$), that
the property of being a real number can be properly defined
in our predicate setting. For short, let us call this theory `FMR' (for `Frege's (theory of) magnitudes (and) real (numbers)'). 
If such an impredicative theory were to be avoided, definition (\ref{5.1}) should be omitted. At most, one could recur 
to a new abbreviation stipulation as
\begin{equation}
\mathsf{R}\left(x\right) := \exists _{\left( \mathcal{P}%
\right) }\mathscr X\exists _{\left( \eth \left( \mathscr X\right) \right)
}R\exists _{\mathscr X}S\left[ x=\mathfrak{R}\left[ R,S\right] \right],  \tag{5.2$'$}
\end{equation}
by then admitting that a real number is an object that satisfies the open formula `$\mathsf{R}\left(x\right)$'.
Call `FMR$'$' the theory got from L$_2$PCA$^{2}_{\Delta^{1}_{0}}$ by merely adding EP$^*$, and replacing 
(\ref{5.1}) with (5.2$'$). The same argument used above to prefer (\ref{3.9}$'$) over (\ref{3.9}) does not apply here, since the definition of real numbers is the final purpose to be reached to revive Frege's program. Hence, if one considers that
this aim is reached only if a property, counting as the property of being a real number, is directly expressed as such,
in the relevant formal setting, on pain of missing the definition itself, there is no other option than going for FMR.

According both to (\ref{5.1}) and (5.2$'$), a real number is a ratio over some
domain of magnitudes. This might appear odd at first glance, since, given different such domains, this might seem to 
entail
that different sorts of real numbers arise,
according to the domain of magnitudes they are defined on. However, from
EP$^*$, it easily follows that, if there are several domains of magnitudes, for
any ratio (or $\mathfrak{R}$-\textit{abstractum}) on one of them, there is
just another ratio (or $\mathfrak{R}$-\textit{abstractum}) on each other of
them that is the same object as the former---i.e. that
the ratio of two first-level binary relations having a certain property $\mathscr M$
such that $\mathcal{M}\left( \mathscr M\right) $ is the same object as the
ratio of two first-level binary relations having another property $\mathscr M^{\prime }$
such that $\mathcal{M}\left( \mathscr M^{\prime }\right) $.

Hence, once real numbers are defined, either in FMR or in FMR$'$, as ratios of magnitudes, 
one can define the usual properties, relations and functions on them,  making the
development of real analysis possible, within these systems---or some appropriate
extensions of them, if needed. We stop here, however, and rather
tackle some meta-theoretical issues concerning
these systems, and the corresponding definitions.

\section{Existence Proofs\label{existence}}

It is easy to see that EP$^*$ implicitly defines continuously many
objects to be identified, either through (\ref{5.1}) or through (5.2$'$), with the real numbers, only in the presence of
nonempty positive classes. If there were no first-level binary relations $\mathbf R$
such that $\exists \mathscr X \left[\mathcal P \left(\mathscr X\right) \wedge \mathscr X \mathbf R\right]$, its 
second-order universal
quantifier would range on an empty domain, and this would
render the right-hand side of EP$^*$ nonsensical, as well as, then, both  (\ref{5.1}) and (5.2$'$) Still, 
a nonempty positive class exists if and only if
this is so for a nonempty domain of magnitudes. Hence, an existence proof of such a domain (or of a positive class)
is an indispensable
supplement to our definition of real numbers: it is required 
to make it sensible. 

Of course, no form of comprehension might be appealed to in order 
to deliver this proof, since no comprehension axiom can secure
the existence of an $R$ such that $\exists \mathscr X \left[ \mathcal M \left(\mathscr X\right) \wedge \mathscr X R\right]$.
Moreover, it would not be enough to observe that the empty first-level binary relation exists no matter what the first-order domain looks like, since, as observed in footnote (\ref{CVDM}), this
relation neither forms nor belongs to a positival (and, then, positive) class. What is to be proved, then, is that there are
enough appropriate (or appropriately related) objects for defining on them continuously many 
(extensionally) distinct first-level binary relations forming a nonempty domain of 
magnitudes\footnote{This is just what Frege seems to signal at the beginning of \S~II.164, in the passage
we have quickly referred to in footnote (\ref{CVDM}) above (Frege 1893-1903, \S~II.164; Frege 2013, pp.~160$_{2}$-61$_{2}$; notice that the English term `Relation' with capital `R' is used here to translate the German term `Relation', which Frege uses, as opposed to `Beziehung', translated instead as `relation', to name value-ranges of relations): 
\begin{quote}
We can now approach the question raised earlier (\S~159) from where we obtain the magnitudes 
whose ratios are real numbers. They will be Relations; and
they must not be empty, i.e., they must not be extensions of relations in which no
objects stand to each other. For such relations have the same extension; there is only
one empty Relation. We could not define any real number with it. If $q$ is the empty
Relation, then [\ldots][$q^{-}$] is the same; likewise [\ldots][$q \sqcup q^{-}$].
Also the composition of Relations on our domain of magnitudes cannot result in the empty Relation; 
but that would happen if there were no object to which some object stood in the first Relation and
which also stood to some object in the second Relation.
We thus require a class of objects that stand to each other in the Relations of our
domain of magnitudes, and, in particular, this class has to contain infinitely many
objects.
\end{quote}}. 
 
This cannot be accomplished by a proof following a similar pattern as the one
that allows neologicists to prove the existence of natural numbers within  the very theory in which they define them, namely 
FA. This proof goes as follows:
\begin{itemize}
\item[--] The concept $\left[ x:x\neq x\right] $ exists by predicative
comprehension;
\item[--] Then, HP allows to define $0$ as $\# \left[x: x\neq x\right] $;
\item[--] By logic, $\left[x: x\neq x\right] \approx \left[x: x\neq x\right]$;
\item[--] Hence, by HP, $0=0$, from which it follows that $0$ 
exists\footnote{Notice that, since HP licenses the
formation of the term `$\# \left[x: x\neq x\right] $', the identity
`$0=0$' might be derived, in classical logic, as an immediate consequence of
the theorem `$\forall x\left[x=x\right]$'. Still, if such a proof of the existence of $0$ were admitted, 
HP would inevitably be endowed with an
existential import that would be incompatible with its alleged analyticity. This is
one of the reasons why it is often advanced that the subjacent logic to FA should be free: the matter has been firstly tackled in Shapiro \& Weir (2000), \S\S~IV-V; but see also Payne (2013).
\label{FL}};
\item[--] Since HP allows to define the successor relation on the whole
first-order domain, natural numbers can be defined as the objects that bear the weak ancestral of this relation to $0$;
\item[--] Proving---from HP plus (impredicative)
comprehension---that any such object has a (unique)
successor is, then, enough to prove, by countable induction, that all the
natural numbers exist.
\end{itemize}

This pattern only allows to prove that there are objects falling under a
first-level concept, given both a way to identify these objects
collectively, as values of a particular function such as $\#$, and a way to
identify them individually, as values of this function for particular
concepts as arguments. In our case, one should, instead, prove that there are
enough objects on which  one can define binary relations
falling under some
second-level concept complying with a certain structural condition, where no particular way is given
to identify both these objects and these relations either 
individually
(except for the identity relation) or collectively. Moreover, by appealing to 
countable induction, this pattern can, at best, be suitable 
for proving the existence of countably many objects, and---even 
if it were
possible to show that such objects allow to define on them
the required binary relations\footnote{We shall hark back on this matter in \S~\ref{Sect.KS} below.}---the
main task of the proof
would just be to prove that, which
is certainly not something that might be done by following this pattern. Hence, though required for making the very definition
of real numbers sensible, the existence proof of nonempty domains of
magnitudes cannot be carried out in the theory FMR itself, and, \textit{a fortiori}, in FMR$'$, in which that definition is stated.

Two alternative
strategies seem possible to deliver it. The first is in line with Frege's 
perspective and looks for an
alternative way to prove the existence of continuously many objects on which
continuously many permutations, forming a domain of
magnitudes, can be defined. The second departs from this perspective, and uses mathematical results unavailable to 
Frege. It might be appealed to, as a sort of unhoped lifeline for Frege's purpose, in order to avoid 
some problems the former strategy suffers from. It consists
of inquiring whether continuously many
permutations forming such a domain can be obtained from
countably many objects, whose existence might, if needed, be proved by applying the
previous proof-pattern. Let us call the first strategy `inflationary' and the
second `non-inflationary'.

\subsection{The Inflationary Strategy: Bicimal Pairs\label{Sect.BP}}
The inflationary strategy can be implemented in at least two slightly different ways, in our setting.
One follows Frege's own plan sketched in  \S ~II.164, and takes the existence of natural numbers for granted. 
The other 
appeals, instead, to a restricted version of BLV, to get an $\omega$-sequence of objects other than Frege's natural numbers.
The structural similarity of these approaches makes them suffer from the same difficulties.
We merely consider the former. The reader might get an idea of the latter from the way 
we deal with a restricted version of BLV at the beginning of \S~\ref{Sect.KS}.

Taking the existence of natural numbers for granted, Frege considers the pairs $<n,\left\{ m_{i}\right\} _{i=0}^{\infty }>$ composed by a natural number and an infinite sequence of positive natural numbers.
These pairs are apt to code Cauchy
series like 
\begin{equation}
n+\sum\limits_{j=1}^{\infty }\lambda _{j}\frac{1}{2^{j}}\hspace{1cm}(\lambda
_{j}\in \left\{ 0,1\right\} ), \label{6.1}
\end{equation}%
under the condition that $m_{i}$ is the $i$-th value of $j$ such that $%
\lambda _{j}=1$ and the $\lambda _{j}$ are not all $0$ after a certain range.
It follows that, once addition is appropriately defined on these pairs, one can associate
to each of them, let us say $\alpha $, a binary (first-level) relation $%
\mathbf{R}_{\alpha }$ such that, for any pair of these same pairs $x$ 
and $y$, $x \mathbf{R}_{\alpha }y $
if and only if $x +\alpha = y$. It is easy to see that this allows to define as many relations as 
pairs like $<n,\left\{ m_{i}\right\} _{i=0}^{\infty }>$ , namely continuously many ones, and that  
these relations are such that:
\begin{itemize}
\item[--] both they and their inverses are
functional, since, for any such pairs $x$, $y$ and $z$,
$x +\alpha = y \wedge x + \alpha = z$ and  $
y +\alpha = x \wedge z + \alpha = x$ each
 entails that $y=z$ ;
\item[--] their composition mimics an addition on the pairs they are defined on, since, if $\beta$ is also such a pair, $\mathbf{R}_{\alpha } \sqcup \mathbf{R}_{\beta }$ extensionally coincides with 
$\mathbf{R}_{\alpha + \beta }$, which is proved by observing that, for any such pairs $x$ and $y$,
$x + (\alpha + \beta) = y$ if and only if there is such a pair $z$ such that
$x+ \alpha = z \wedge
z + \beta = y$;
\item[--] the identity relation is not one of them, since no Cauchy series like (\ref{6.1}) is equal to zero. 
\end{itemize}
It is, then, easy to verify that these relations form a positive class, from which a domain of magnitudes is obtained by merely adding their inverses and the identity relation.

Objects rendering these pairs in a formal setting are quite easy to define
in any second-order version of arithmetic. A simple way to do it (Panza, 2016 and Panza, FC3) is by
adding a new axiom, under the form of the following
abstraction principle:
\begin{equation*}
\forall _{\mathbf{N}}XY\forall _{\mathbf{N}}xy\left[ <x,X>\,=\,<y,Y>\,%
\Leftrightarrow \left( x=y\wedge \forall z\left( Xz\Leftrightarrow Yz\right)
\right) \right],\tag{FP}
\end{equation*}%
where the index `$\mathbf{N}$' signals that the universal quantifiers are restricted to properties of natural numbers and to these 
very numbers respectively---the acronym `FP' stands fro  `Frege's Principle', and emphasizes the fact that this principle is 
a restricted adapted form of BLV.

Of course, to go ahead, we have to prove that FP is consistent. Assuming the consistency
of second-order arithmetic, to this purpose, it is, however, enough to
observe that FP has a model in the V$_{\omega +1}$ stage of ZF's hierarchy. This is because second-order arithmetic has a model in the V$_{\omega}$ segment of ZF, and
consequently the set $\wp\omega$ of all subsets of the set of natural numbers is at stage V$_{\omega +1}$, and 
provides the required model.

This having been established, we can look at the pairs like $<n,\mathbf P>$, implicitly
defined by FP, and formed by a natural number $n$ and a property $\mathbf P$ of natural numbers, and select among them
those whose second element $\mathbf P$ is an
infinite property of natural numbers---i.e. it is such that
$\forall _{\mathbf{N}}x\left[ \mathbf P x\Rightarrow \exists _{\mathbf{N}}y\left[
x<y\wedge \mathbf P y\right] \right] $. For short, call them `bicimal pairs'. Clearly, FP allows to distinguish continuously
many such pairs. 

They can be arranged into two partitions, such that any bicimal pair $<n,\mathbf P>$
belongs to one partition if $\mathbf P0$, and to the other if $\lnot \mathbf P0$. A
total order can, moreover, be defined on these pairs, in such a way that the pairs in 
the former partition count
as positive, the pair  $<0,\mathbf N^{+}>$ (where `$\mathbf N^{+}$' designates the property of being a positive natural number)
counts as the zero pair, and the other pairs in the latter partition count as negative (more details are given in Panza, 2016, p.~417; others
will be found in Panza, FC3). This makes the bicimal pairs form an
additive Abelian group, that can be proved to be dense,
totally ordered and Dedekind-complete (and, then, Archimedean), and can also be extended to a
field by an appropriate definition of multiplication (details are, again, to be found in Panza, FC3).
It would, then, be not only very tempting, but also rather natural to
code the real numbers by bicimal 
pairs, so as to avoid the very
definition of domains of magnitudes and of ratios thereof as perfectly useless. 

Still, this is certainly not what Frege's strategy should lead us to. In order to follow his indications, one should rather 
define appropriate permutations on bicimal pairs and show that they form a domain of magnitudes. This can easily be done by associating
to  any such pair $<n,\mathbf P>$ the binary relation $\mathbf R_{<n,\mathbf P>}$ such that, for any
two such pairs $<y,Y>$ and $<z,Z>$, $<y,Y>\mathbf R_{<n,\mathbf P>}<z,Z>$ if and only if $%
{<y,Y>}+{<n,\mathbf P>} = {<z,Z>}$. It would, then, be easy to verify that the relations associated with
positive bicimal pairs just behave as those Frege suggests to associate to his pairs,
and form, then, a positive class, which is easy to extend to a domain of magnitudes.

Following this path leads, then, to an arithmetical copy of the additive (ordered) group of the real
numbers, as an intermediate step in a much more complex, supposedly
non-arithmetical definition of them. Hence, real numbers might be \textit{ipso facto} identified
with bicimal pairs, by so dramatically departing, however, from Frege's purpose. 
The same happens for any other way of pursuing the inflationary strategy: it
inevitably leads either to encode the real numbers by objects other than ratios of magnitudes and thus depart from Frege, 
or to define real numbers twice, structurally speaking, by
accepting the idea that the second definition requires a supplementary axiom which is not at 
all required by the first, namely  EP$^*$.

Though mathematically
quite shocking, the former option poses no
problem from a realist perspective such as Frege's,
since the realist may argue that bicimal pairs are intrinsically not real numbers, though they behave like them.
We do not want to dig into this possibly rather odd attitude. We merely observe that, in this perspective, the ratios on the appropriate permutations defined on these pairs 
could not but be taken to be real numbers. Thus, 
the only way to avoid concluding
that, \textit{pace} Frege, real numbers are intrinsically arithmetical objects would be to prove that
there are non-arithmetical nonempty domains of magnitudes. Insofar as  EP$^*$ identifies
the ratios on any domain
of magnitudes
with the ratios on any other such domain, this would leave room for arguing that being a real number
is not intrinsically the same as being a ratio on permutations defined on bicimal pairs, since
such a ratio would merely provide one among other possible and
essentially distinct
modes of presentation of such a number.
But, then, how to prove the existence of other, non-arithmetical nonempty domains of magnitudes?

\subsection{The Non-Inflationary Strategy\label{Sect.KS}}

As\footnote{%
Both for this \S\ and the following one, we are very much
indebted to Mirna D\v{z}amonja, Andrew Moshier and, overall, Alain
Genestier who guided us in the understanding of Karrass and Solitar's proof
and annexed topics.} a matter of fact, also the non-inflationary strategy might be grounded on the assumption of the existence of 
natural numbers. Strictly speaking, this is not necessary, however: it might also be grounded on
a consistently restricted version of BLV.

Let `$\mathscr F(\mathbf{X})$' be the
abbreviation of a logical
second-order formula stating that $\mathbf X$ is
a property satisfied at most by finitely many objects. A
possibility is appealing to Dedekind-finiteness:
\[
\mathscr F(X) := \forall Y \left[
\begin{array}{l}
\forall x \left[
\left[Yx \Rightarrow Xx \right] \wedge \exists y \left[Xy \wedge \lnot Yy
\right]
\right] \Rightarrow \medskip \\
\lnot \exists R \left[ 
\begin{array}{l}
\forall z \left[Xz \Rightarrow \exists ! z' \left(
zRz' \wedge Yz'
\right)
\right] \wedge \medskip \\
\forall z' \left[Yz' \Rightarrow \exists ! z \left(
zRz' \wedge Xz
\right)
\right]
\end{array}
\right]
\end{array}
\right].
\]
The
relevant restricted version of BLV, call it `FinBLV', is, then, this:
\begin{equation*}
\forall _{\left( \mathscr F\right) }XY\left[ \varepsilon
X=\varepsilon Y\Leftrightarrow \forall x\left( Xx\Leftrightarrow Yx\right) %
\right], \tag{FinBLV}
\end{equation*}
where `$\varepsilon X$' and `$\varepsilon Y$' denote the
extensions of $X$ and $Y$, and the index `$\mathscr F$' signals that the universal
quantifier is restricted to the (first-level) properties
satisfying the formula `$\mathscr F(X)$'\footnote{This principle is freely inspired by two different suggestions respectively advanced by
Antonelli \& May (2005, p.~11), and by Boolos (1998a, p.~99; 1998b, p.~178)
, in particular by Boolos' \textit{New V}, i.e. $\forall F\forall G(\varepsilon F=\varepsilon G\leftrightarrow (Small(F) \vee Small(G)\rightarrow F\equiv G))$, where `\textit{Small}' means `not equinumerous with the universal concept $[x: x=x]$'. By remaining faithful to this suggestion, we should replace `\textit{Small}' by `Fin' (or take the former to mean `finite'), and FinBLV by ` $\forall F\forall G(\varepsilon F=\varepsilon G\leftrightarrow (\mathscr F (F) \vee \mathscr F (G)\rightarrow F\equiv G))$'. Though this latter principle would not be equivalent to FinBLV, we cannot see any
relevant difference between them with respect to our purpose. We prefer FinBLV simply because it makes immediately clear that infinite concepts do not matter, here, by ascribing to them no extension, rather
than ascribing to all of them the same extension.}.

FinBLV implicitly defines a countable infinity of extensions. To see it, notice that $[x:x\neq x]$ is a finite property, whose existence is warranted by
predicative second-order comprehension. So FinBLV applies to it. Let $%
\oslash $ be its extension. Since it is a theorem of second-order (free)\footnote{See footnote (\ref{FL}), above.} logic
that $\forall x([x:x\neq x]\Leftrightarrow \lbrack x:x\neq x])$, from FinBLV
it follows that $\oslash =\oslash $, which entails
that $%
\exists x\left( x=\oslash \right) $. 
One can, then, firstly appeal to second-order impredicative comprehension to
define a functional first-level binary relation on the whole
first-order domain, 
by stating that 
\[
\forall x,y\left[x\mathsf{S}y\Leftrightarrow
\lbrack \exists X_{\mathscr F}\left[ x=\varepsilon X\right] \wedge y=\varepsilon \left[
z:z=x\right]\right],
\]
then define the weak ancestral of $\mathsf{S}$, i.e. $%
\mathsf{S}^{\ast =}$, and finally appeal to second-order predicative comprehension to define
the property $\mathsf{E}_{\oslash }$ of being an extension
belonging to
the $\mathsf S$-succession starting from $\oslash $:
\[
\forall x[%
\mathsf{E}_{\oslash }x\Leftrightarrow \oslash \mathsf{S}^{\ast =}x].
\]
This allows to accomplish the task by repeating, \textit{mutatis mutandis}, the neologicist recursive proof
of the existence of natural numbers\footnote{The proof depends on the lemma that $\forall x\exists y\left[
y=\varepsilon [z:z=x]\right] $. Here is how this can be proved. FinBLV and $\forall X\forall x \left[Xx \Leftrightarrow Xx \right]$ 
imply that $\forall_{\left( \mathscr F\right) } X\left[\varepsilon X = \varepsilon X\right]$, and, then,
$\forall X_{\left( \mathscr F\right) } \exists y \left[y= \varepsilon X \right]$. What is required is, then, 
to prove that `$[z : z = x]$' denotes a (finite) property, which is ensured by predicative
comprehension with parameters, since it entails that $\forall x\exists X\forall z\left[
Xz\Leftrightarrow z=x]\right] $. Notice that this
proof also holds in (negative) free logic: thanks to Ludovica Conti for drawing our attention to this; see also Conti (2019, p. 145, fn. 426, and pp. 151-152).
This lemma implies, \textit{a fortiori%
}, that $\forall _{\mathsf{E}_{\oslash }}x\exists y\left[ y=\varepsilon
\lbrack z:z=x]\right] $. The very definition of the weak ancestral of $%
\mathsf{S}$ allows, then, to prove quite easily the principle of induction
for the FinBLV-\textit{abstracta} having the property $\mathsf{E}_{\oslash }$ (or
$\mathsf{E}_{\oslash }$-\textit{abstracta}, from now on)---namely 
`$\forall X\left[
\left( X\oslash \wedge X\underset{\mathsf{E}_{\oslash }}{\overset{%
\mathcal{S}}{\longrightarrow }}X\right)
\Rightarrow \forall _{\mathsf{E}_{\oslash }}x\left[ Xx\right] \right] $', where `$X\underset{\mathsf{E}_{\oslash }}{\overset{\mathsf{S}}{\longrightarrow }}X$' abbreviates `$\forall_{\mathsf{E}_{\oslash }} x\forall y%
\left[ \left( Xx\wedge x\mathsf{S}y\right) \Rightarrow Xy\right] $'.\ With
this principle at hand, it is, then, equally easy to prove that $\forall _{%
\mathsf{E}_{\oslash }}x\exists X\left[ x=\varepsilon X\right] $, as a
consequence of the stipulation that $\oslash =\varepsilon \lbrack x:x\neq x]$%
, from which it immediately follows that $\exists X\left[ \oslash
=\varepsilon X\right] $. By appealing to $\forall _{\mathsf{E}_{\oslash }}x\exists y\left[ y=\varepsilon
\lbrack z:z=x]\right] $, one gets that $%
\forall _{\mathsf{E}_{\oslash }}x\exists y\left[ \exists X\left[
x=\varepsilon X\right] \wedge y=\varepsilon \lbrack z:z=x]\right] $, that
is, $\forall _{\mathsf{E}_{\oslash }}x\exists y\left[ x\mathsf{S}y\right] $.
Next, the properties of the ancestral of a binary relation allow to prove
that $\forall _{\mathsf{E}_{\oslash }}x\forall y\left[ x\mathsf{S}%
y\Rightarrow \mathsf{E}_{\oslash }y\right] $, and so to conclude that $%
\forall _{\mathsf{E}_{\oslash }}x\exists _{\mathsf{E}_{\oslash }}y\left[ x%
\mathsf{S}y\right] $. The existence of countably many 
$\mathsf{E}_{\oslash }$-\textit{abstracta} finally follows from
proving by induction that $\forall _{\mathsf{E}_{\oslash
}}x[\lnot x\mathsf{S}^{\ast \neq }x]$, where $\mathsf{S}^{\ast \neq }$ is
the strong ancestral of $\mathsf{S}$. This last part of the proof is a
little bit harder than the previous ones, but requires no specific skills:
it just parallels the analogous proof of FA, by exploiting, like this
(together with the principle of induction and some properties of the strong
ancestral, also) the obvious facts that $\lnot \oslash \mathsf{S}^{\ast \neq
}\oslash $ and $\forall x,y,z\left[ \left( x\mathsf{S}y\wedge z\mathsf{S}%
y\right) \Rightarrow x=z\right] $. Though quite simple, this proof
of existence of countably many
$\mathsf{E}_{\oslash }$-\textit{abstracta} directly
involves most of Peano's second-order axioms as theorems about them. The axioms that do not enter it, i.e. `$\lnot \exists _{\mathsf{E}%
_{\oslash }}x\left[ 0\mathsf{S}x\right] $' and `$\forall x,y,z\left[ \left( x%
\mathsf{S}y\wedge x\mathsf{S}z\right) \Rightarrow y=z\right] $' (the second
of which would allow replacing the relation $\mathsf{S}$ by the successor function), are, moreover, even easier to prove. Hence, if
impredicative comprehension is accepted, Peano's second-order arithmetic is
interpretable on the mentioned extensions---with no appeal to set-theoretical
assumptions on them. These extensions do not comply,
however, with HP, which makes them crucially differ  from
natural numbers as defined in FA.\label{ftn:test}}.

Now, consider the symmetric group $\mathbf{\Sigma }_{\mathbb{N}}$ on
the natural numbers, i.e. the (additive) group of all permutations on $%
\mathbb{N}$. We know that $\mathbf{\Sigma }_{\mathbb{N}}$ has cardinality $
2^{\aleph _{0}}$. But we also know that it contains torsion elements%
\footnote{%
A torsion element of a group $\mathbf{G}$ is an element $g$ of $\mathbf{G}$
such that $g^{n}=e$ for some natural number $n$, where $e$ is the identity
element of $\mathbf{G}$.}, which prevents both from defining on $\mathbf{%
\Sigma }_{\mathbb{N}}$ a total order compatible with the group structure,
and from making an injective map from a (nonempty) domain of magnitudes (if
any) into it surjective. It follows that $\mathbf{\Sigma }_{\mathbb{N}}$
does not provide a nonempty model for our Fregean consistent definition of
such a domain. Nevertheless, 
by a theorem by
Karrass and Solitar (1956, p.~65), $\mathbf{%
\Sigma }_{\mathbb{N}}$ provably \textquotedblleft contains a copy of the
additive group of the reals\textquotedblright . In other terms, there is a
subgroup of $\mathbf{\Sigma }_{\mathbb{N}}$ which is isomorphic to $\left( 
\mathbb{R},+\right) $, and is, then, not only totally ordered, but also
Abelian, dense, Archimedean and Dedekind-complete. 
Since any totally ordered, dense and Dedekind-complete group of
permutations is a model of our definition of domains of magnitudes\footnote{More precisely, 
the second-order property of being a permutation
belonging to any such group satisfies the right-hand
side of (3.6$'$).}, $\mathbf{\Sigma }_{\mathbb{N}}$ contains such a model. Insofar as $\mathbf{\Sigma 
}_{\mathbb{N}}$ is isomorphic to the symmetric group on whatever infinite countable set,
it is also so to the symmetric group $\mathbf{\Sigma 
}_{\mathsf{E}_{\oslash }}$
formed by the $\mathsf{E}_{\oslash }$-\textit{abstracta}. Hence, 
Karras and Solitar's theorem entails that admitting the
existence of this latter symmetric group ensures the existence of
a nonempty domain of magnitudes. It
would, then, be enough to admit that it makes sense to speak of
all permutations on an infinite countable domain $\mathbf{D}$ of objects\footnote{As Frege himself should have admitted, 
in order to make his own
definition of domains of magnitudes sensible.}, and that
the existence of (the elements of) $\mathbf{D}$
\textit{ipso facto} entail the existence of the group formed by these
permutations,
to conclude that the existence of a
nonempty domain of magnitudes is ensured by the existence of the natural numbers
or of the $\mathsf{E}_{\oslash }$-\textit{abstracta}, because of
an immediate corollary of Karrass and Solitar's result\footnote{Notice that though the definition of $\mathsf{S}^{\ast =}$ allows proving that the
$\mathsf{E}_{\oslash }$-\textit{abstracta} form an $\omega$-sequence, as showed in footnote (\ref{ftn:test}) above,
all that is relevant here is the cardinality of set formed by these \text{abstracta}, namely the fact that this set is countably infinite.}.

Karrass and Solitar's proof could certainly not have been within Frege's reach. But it is not very entangled, as such, and, under
the mentioned admission, most of it can be conducted constructively, which is something
Frege seems to require for his proof.

Let $\mathbf{I}$ be an infinite countable set, for example an infinite
subset of $\mathbb{N}$. Let $\mathbf{\sigma }=\coprod\limits_{i\in \mathbf{I}%
}\mathbf{C}_{i}$ be a partition of $\mathbb{N}$ in infinite (countable)
subsets $\mathbf{C}_{i}$ ($i\in $ $\mathbf{I}$), which makes $%
\bigcup\limits_{i\in \mathbf{I}}\mathbf{C}_{i}=\mathbb{N}$ and $\mathbf{C}%
_{h}\cap \mathbf{C}_{k}=\mathbb{\emptyset }$, for any distinct $h,k$ in $%
\mathbf{I}$. Let $\mathbf{\varrho }:\mathbb{N}\longrightarrow \mathbf{I}$ be
the (surjective) application defined by this partition, associating to any $n
$ in $\mathbb{N}$, the single element $i=\mathbf{\varrho }\left( n\right) $
of $\mathbf{I}$ such that $n\in \mathbf{C}_{i}$. Such a partition, and
therefore such an application, can easily be defined constructively.\ To
make a simple example, take $\mathbf{I}$ to be the set of all prime numbers
plus $0$, $\mathbf{C}_{0}$ the set of all natural numbers that are not
(positive) powers of a prime number, namely $\left\{ 0,1,6,10,\ldots
\right\} $, and, for any prime number $p$, $\mathbf{C}_{p}$ the set of all
(positive) powers of $p$. Though we would not be (presently) able to write a
general formula providing, for any natural number $n$, the value $\mathbf{%
\varrho }\left( n\right) $ of $i$, such that $n\in S_{\mathbf{\varrho }%
\left( n\right) }$, there is a finite algorithm allowing us to calculate
such a value $\mathbf{\varrho }\left( \mathbf{n}\right) $ for whatever given
natural number $\mathbf{n}$.

For any $i$ in $\mathbf{I}$, let $\mathbf{\pi }_{i}$ be a permutation on $%
\mathbf{C}_{i}$. Define $\mathbf{\pi }:\mathbb{N}\longrightarrow \mathbb{N}$
by establishing that for any $n$ in $\mathbb{N}$, $\mathbf{\pi }\left(
n\right) =\mathbf{\pi }_{\mathbf{\varrho }\left( n\right) }\left( n\right) $%
. This is clearly a permutation on $\mathbb{N}$, so that $\mathbf{\pi }\in 
\mathbf{\Sigma }_{\mathbb{N}}$. If $\coprod\limits_{i\in \mathbf{I}}\mathbf{S%
}_{i}$ and $\mathbf{\varrho }$ have been defined constructively,  any $%
\mathbf{\pi }_{i}$, and therefore $\mathbf{\pi }$, can also be so defined.
Supposing $\mathbf{\varrho }$ be defined as in the previous example, we
might, for example, take $\mathbf{\pi }_{i}$ to be the permutation
exchanging the $\left( 2j-1\right) $-th element of any set $\mathbf{C}_{i}$,
according to the usual order on $\mathbb{N}$, with the $2j$-th one, and vice
versa ($j=1,2,\ldots $). Then $\mathbf{\pi }$ would permute any natural
number with another natural number following or preceding it in this subset,
according to this order: $%
\mathbf{\pi }\left( 0\right) =1$, $\mathbf{\pi }\left( 1\right) =0$, $%
\mathbf{\pi }\left( 2\right) =4$, \ldots\ Insofar as the same can be done
for any permutation $\pi _{i}$ on any $\mathbf{C}_{i}$, the application $\pi
:\mathbb{N}\longrightarrow \mathbb{N}$ defined by stating that $\pi \left(
n\right) =\pi _{\mathbf{\varrho }\left( n\right) }\left( n\right) $ provides
a group monomorphism\footnote{%
A group monomorphism is an injective group homomorphism, i.e. an injective
map $\mu $ from a group $\left( \mathbf{G},\ast \right) $ to another group $%
\left( \mathbf{H},\star \right) $, such that $\mu \left( x\ast y\right) =\mu
\left( x\right) \star \mu \left( x\right) $ for any $x$, $y$ in $\mathbf{G}$.}%
 $\prod\limits_{i\in \mathbf{I}}\mathbf{\Sigma }_{\mathbf{C}%
_{i}}\longrightarrow \mathbf{\Sigma }_{\mathbb{N}}$, where $%
\prod\limits_{i\in \mathbf{I}}\mathbf{\Sigma }_{\mathbf{C}_{i}}$ is the
product of the symmetric groups on all sets $\mathbf{C}_{i}$---since, if $%
\left\{ \pi _{i}\right\} _{i\in \mathbf{I}}$ and $\left\{ \pi _{i}^{\prime
}\right\} _{i\in \mathbf{I}}$ are two families of permutations on all the
sets $\mathbf{C}_{i}$, then $\left\{ \pi _{i}\right\} _{i\in \mathbf{I}%
}\circ \left\{ \pi _{i}^{\prime }\right\} _{i\in \mathbf{I}}=\left\{ \pi
_{i}\circ \pi _{i}^{\prime }\right\} _{i\in \mathbf{I}}$. Under the
condition mentioned above, and provided all sets $\mathbf{C}_{i}$ be
(constructively) given, this further step of the proof is, also,
constructively licensed.

By Cayley's theorem, any group $\mathbf{G}$ is isomorphic to a subgroup of
the symmetric group $\mathbf{\Sigma }_{\mathbf{G}}$ on $\mathbf{G}$ itself.
Thus, there is 
a group monomorphism $\left( \mathbb{Q},+\right)
\longrightarrow \mathbf{\Sigma }_{\mathbb{Q}}$ from the additive group of
the rational numbers $\left( \mathbb{Q},+\right) $ into the symmetric group $%
\mathbf{\Sigma }_{\mathbb{Q}}$ on $\mathbb{Q}$. Though quite general,
Cayley's theorem can easily be proved constructively, which also makes this new step of
Karrass and Solitar's proof admissible from Frege's perspective. A quite usual way to prove it 
is, for example, by considering the translation $\tau
_{y}:x\longmapsto y \ast x$ on the domain of $\mathbf{G}$ (where $ \ast$ is the
composition law of this group), since it is
easy to see that $\tau _{(a \ast b)}=\tau _{a}\circ \tau _{b}$, for any $a$, $b$
in such a domain. This proof directly applies to the present case, for $\mathbf{G}$
is nothing but $\left( \mathbb{Q},+\right) $. To this purpose, we can take  
$\mathbb{Q}$ to play the role of the domain of $\mathbf{G}$ and $+$ that of
$\ast $, and observe that $\tau _{(q+r)}=\tau _{q}\circ \tau _{r}$, for any $q$, $r$
in $\mathbb{Q}$. Notice, moreover, that this application is perfectly akin
to Frege's suggested definition of permutations on the pairs $<n,\left\{
m_{i}\right\} _{i=0}^{\infty }>$ in the outline of his
existence proof.\footnote{%
See \S ~\ref{Sect.BP}, above.}

For any $i$ in $\mathbf{I}$, let us choose a bijection $\mathbf{\vartheta }_{i}:%
\mathbb{Q}\longrightarrow \mathbf{C}_{i}$ from the set $\mathbb{Q}$ to the
set $\mathbf{C}_{i}$.\ Because of the monomorphism $\left( \mathbb{Q}%
,+\right) \longrightarrow \mathbf{\Sigma }_{\mathbb{Q}}$, this engenders,
for any such $i$, a new group monomorphism $\left( \mathbb{Q},+\right)
\longrightarrow \mathbf{\Sigma }_{\mathbf{C}_{i}}$ from $\left( \mathbb{Q}%
,+\right) $ into the symmetric group $\mathbf{\Sigma }_{\mathbf{C}_{i}}$ on
any $\mathbf{C}_{i}$, and, by composition, a further group monomorphism $%
\prod\limits_{i\in \mathbf{I}}\left( \mathbb{Q},+\right) \longrightarrow
\prod\limits_{i\in \mathbf{I}}\mathbf{\Sigma }_{\mathbf{C}_{i}}$ from the
product $\prod\limits_{i\in \mathbf{I}}\left( \mathbb{Q},+\right) $ of
countably many copies of the additive group $\left( \mathbb{Q},+\right) $
into the product $\prod\limits_{i\in \mathbf{I}}\mathbf{\Sigma }_{\mathbf{C}%
_{i}}$. Again, if the partition $\mathbf{\pi }$ and the application $\mathbf{%
\varrho }$ are defined constructively, making the choice of the
bijections $\mathbf{\vartheta }_{i}$ and so getting these two monomorphisms
require no form of the Axiom of Choice, and, under the condition mentioned
above, is unquestionably constructive.
By combining the
monomorphisms $\prod\limits_{i\in \mathbf{I}}\mathbf{\Sigma }_{\mathbf{C}%
_{i}}\longrightarrow \mathbf{\Sigma }_{\mathbb{N}}$ and $\prod\limits_{i\in 
\mathbf{I}}\left( \mathbb{Q},+\right) \longrightarrow \prod\limits_{i\in 
\mathbf{I}}\mathbf{\Sigma }_{\mathbf{C}_{i}}$, we finally get,
constructively again, a final monomorphism%
\begin{equation}
\prod\limits_{i\in \mathbf{I}}\left( \mathbb{Q},+\right) \longrightarrow 
\mathbf{\Sigma }_{\mathbb{N}}.  \label{GMQR}
\end{equation}

This provides a constructive ground for Karrass and Solitar's proof. But, for its completion, 
a last step is needed, which,
instead, requires non-constructive means and makes, then, the whole proof non-constructive.
Both additive groups $%
\left( \mathbb{Q},+\right) $ and $\left( \mathbb{R},+\right) $, enriched
with the usual multiplication by a rational number, can be regarded as
vector spaces over the field $\left( \mathbb{Q},+,\cdot \right) $, and this
is also the case for the product $\prod\limits_{i\in \mathbf{I}}\left( 
\mathbb{Q},+\right) $. When $\prod\limits_{i\in \mathbf{I}}\left( \mathbb{Q}%
,+\right) $ and $\left( \mathbb{R},+\right) $ are so regarded, it is
however not plain that they have a basis, unless Zorn's lemma is
appealed to, since this lemma allows to prove that every vectorial
space has a basis\footnote{%
By speaking of basis of a vector space, we more precisely mean, here, a
Hamel basis. Let $\mathbf{V}$ be a vector space on a field $\mathbf{F}$. An
Hamel basis of $\mathbf{V}$ is a set $\mathbf{B}_{\mathbf{V}}$ of linearly
independent vectors in $\mathbf{V}$, such that for any vector $v$ in $%
\mathbf{V}$ there is a (unique) finite subset $\left\{ \mathbf{v}_{1},%
\mathbf{v}_{2},\ldots \mathbf{v}_{k}\right\} $ of $\mathbf{B}_{\mathbf{V}}$
such that $v=\mathbf{a}_{1}\mathbf{v}_{1}+\mathbf{a}_{2}\mathbf{v}%
_{2}+\ldots +\mathbf{a}_{k}\mathbf{v}_{k}$, where $\mathbf{a}_{1}$, $\mathbf{%
a}_{2}$\textbf{, \ldots , }$\mathbf{a}_{k}$ are elements of $\mathbf{F}$.}.
The non-constructive step of the proof just consists, then, in appealing
to this lemma
for proving that  $\prod\limits_{i\in 
\mathbf{I}}\left( \mathbb{Q},+\right) $ and $\left( \mathbb{R},+\right) $
have a basis. 

This makes it possible to appeal to a new theorem ensuring that, if
a vector space has several distinct
bases, all of them have the same cardinality---which is, then, to be regarded
as the dimension of this space. For vector spaces with
finite bases, this is quite easy to prove\footnote{A simple combinatorial argument allows to prove that the cardinality 
of any set of
linearly independent vectors is smaller or equal to the cardinality of
whatever generating set. Insofar 
as a basis of a vector space is a set of
linearly independent vectors that generates the whole space, from this it
immediately follows that two bases cannot have different
cardinality, since, if they did, there would be a set of linearly
independent vectors whose cardinality is greater than that of a generating
set.}. For vector spaces whose generating sets are infinite, as 
$\prod\limits_{i\in \mathbf{I}}\left( \mathbb{Q}%
,+\right) $ and $\left( \mathbb{R},+\right) $, the proof is more entangled, but 
can still be given constructively. In this case, the theorem can indeed be viewed as 
an immediate corollary
of another theorem asserting that the cardinality of any generating set of
a vector space $\mathbf{V}$ that can be regarded as
the direct sum of an infinite family $\left\{\mathbf{V}_{\lambda \in \mathbf{\Lambda}}\right\}$ of
non-zero vectorial sub-spaces is greater or equal to the cardinality of the set of indices $\mathbf{\Lambda}$
(Bourbaki, Algebra I, ch. II, prop. 23, cor. 2).

Provided that two vector spaces on the same field (both having bases) are isomorphic if (and only if) they
have the same dimension,
$\prod\limits_{i\in 
\mathbf{I}}\left( \mathbb{Q},+\right) $ and $\left( \mathbb{R},+\right) $, when regarded as
vector spaces over $\left( \mathbb{Q},+,\cdot \right) $,
have the same dimension, and thus 
are isomorphic. This
obviously entails that the group monomorphism (\ref{GMQR}) results in a new
group monomorphism 
\begin{equation*}
\left( \mathbb{R},+\right) \longrightarrow \mathbf{\Sigma }_{\mathbb{N}},
\end{equation*}%
which makes $\mathbf{\Sigma }_{\mathbb{N}}$ contain a copy of $\left( 
\mathbb{R},+\right) $, as was required to be proved.

\label{CVFinBVL}If this proof is granted, together with the existence both of an infinite
countable set $\mathbf{D}$---as $\mathbb{N}$ or the set formed by the $\mathsf{E}_{\oslash }$-\textit{abstracta}---and of the symmetric
group $\mathbf{\Sigma }_{\mathbf{D}}$ on this set, the conclusion follows
that there is (at least) a totally ordered, dense and Dedekind-complete
subgroup of $\mathbf{\Sigma }_{\mathbf{D}}$. Let $\left( \mathbf{M}_{\mathbf{%
D}},\circ ,<\right) $ be such a subgroup of $\mathbf{\Sigma }_{\mathbf{D}}$.
Claiming that $\left( \mathbf{M}_{\mathbf{D}},\circ ,<\right) $ is a domain
of magnitudes in agreement with definition (\ref{3.9}$'$) is the same as
claiming that the triple $\left( \mathscr M_{\mathbf{D}},\sqcup ,\mathscr P_{%
\mathbf{D}}\right) $ satisfies the open formula `$\left[ \mathcal{P}\left( %
\mathscr Y\right) \wedge \forall R\left[ \mathscr XR\Leftrightarrow \eth %
\left( \mathscr Y\right) \left( R\right) \right] \right] $' entering the
right-hand side of this definition, with `$\mathscr M_{\mathbf{D}}$' as a
value of `$\mathscr X$' and `$\mathscr P_{\mathbf{D}}$' as a value of `$%
\mathscr Y$', provided that: any binary relation $R\ $has the property $%
\mathscr M_{\mathbf{D}}$ if and only if it is in $\mathbf{M}_{\mathbf{D}}$,
and the property $\mathscr P_{\mathbf{D}}\ $only if it has the property $%
\mathscr M_{\mathbf{D}}$; `$\sqcup $' denotes the same operation on the
binary relations that are in $\mathbf{M}_{\mathbf{D}}$ as that denoted by `$%
\circ $'; and, for any $R,S$ in $\mathbf{M}_{\mathbf{D}}$, $\mathscr P_{%
\mathbf{D}}R\sqcup S^{-}$ if and only if $S<R$, so that $\mathscr P_{\mathbf{%
D}}R$ if and only if $0_{\mathbf{M}_{\mathbf{D}}}<R$---where $0_{\mathbf{M}_{%
\mathbf{D}}}$ is, of course, the neutral element of $\left( \mathbf{M}_{%
\mathbf{D}},\circ ,<\right) $, namely the identity relation. To make this
claim sensible, we have, of course, to grant that these conditions ensure
the existence of the two second-level properties $\mathscr M_{\mathbf{D}}$
and $\mathscr P_{\mathbf{D}}$, which requires third-order comprehension.
Supposing it is admitted, Karrass and Solitar's result provably
establishes that, under the admission of the existence of $\mathbf{D}$
and $\mathbf{\Sigma }_{\mathbf{D}}$, there is a nonempty domain of magnitudes, namely
the ordered group defined by the triple $\left( \mathscr M_{\mathbf{D}%
},\sqcup ,\mathscr P_{\mathbf{D}}\right) $, since it entails that the
properties $\mathscr P_{\mathbf{D}}$ and $\mathscr M_{\mathbf{D}}$ are such
that $\mathcal{P}\left( \mathscr P_{\mathbf{D}}\right) $ and $\forall R\left[
\mathscr M_{\mathbf{D}}R\Leftrightarrow \eth (\mathscr P_{\mathbf{D}})\left(
R\right) \right] $, so that $\mathcal{M}\left( \mathscr M_{\mathbf{D}%
}\right) $.

This being granted, it is, then, easy to prove that there are as many 
distinct binary relations in such a domain as distinct ratios defined on it according to 
EP$^*$, namely that these ratios are
continuously many. In other terms, the ordered pairs $\left[ R,S\right] $ of binary relations,
the first of which has $\mathscr M_{\mathbf{D}}$ and the second $\mathscr P_{%
\mathbf{D}}$, form continuously many distinct equivalence classes according
to the equivalence relation on the right-hand side of EP$^*$, under the
replacement of both `$\mathscr X$' and `$\mathscr Y$' with `$\mathscr P_{%
\mathbf{D}}$'\footnote{%
Notice that from the conditions above, it immediately follows that $\eth %
\left( \mathscr P_{\mathbf{D}}\right) \left( R\right) $ if and only if $%
\mathscr M_{\mathbf{D}}R$.}. For shortness and simplicity, let us
sketch this proof with reference to EP. Its rephrasing with respect to
EP$^*$ is only a matter of laborious routine. 

Consider the first of the three
disjoints forming the right-hand side of EP.
Suppose that there be three binary relations $\mathbf R$, $\textbf S$, 
and $\textbf S^{\prime }$ that have $\mathscr P_{%
\mathbf{D}}$ and are such that 
\[
\begin{array}{l}
\mathbf S \sqsubset_{\mathscr P_{\mathbf{D}}} \mathbf S^{\prime}, \medskip \\
\forall _{\mathbf{N}}x,y\left[ x\mathbf R\sqsubset_{\mathscr P_{%
\mathbf{D}}}y\mathbf S\Leftrightarrow x\mathbf R\sqsubset_{\mathscr P_{\mathbf{D}%
}}y\mathbf S^{\prime }\right].
\end{array}
\]
If there were, then, another
binary relation $\textbf T$ such that $\mathbf S^{\prime }$ coincided with $\mathbf S \sqcup \mathbf T$, it would follow that 
\[
\forall _{\mathbf{N}}x,y\left[ x\mathbf R \sqsubset_{\mathscr P_{\mathbf{D}}}y \mathbf S\Leftrightarrow x \mathbf R \sqsubset_{%
\mathscr P_{\mathbf{D}}}y\left(\mathbf S\sqcup \mathbf T\right)\right],
\]
 which is impossible
because of the Archimedeanicity of $\left( \mathbf{M}_{\mathbf{D}},\circ
,<\right) $.\ This proves that, for whatever binary relation $\mathbf R$ that has 
$\mathscr P_{\mathbf{D}}$, there are as many distinct ratios 
$\mathfrak{R}\left[ \mathbf R,S\right]$ (where $S$ is
a binary relation that has 
$\mathscr P_{\mathbf{D}}$) as binary relations that have 
$\mathscr P_{\mathbf{D}}$,
 namely continuously many such ratios. 
Consider a binary relation $\mathbf R^{\prime} $ that has 
$\mathscr P_{\mathbf{D}}$ and is distinct from $\mathbf R$. Because of the completeness of positive classes, 
there is one and only one binary relation $S$ that also has  $\mathscr P_{\mathbf{D}}$ such that
\[
\forall _{\mathbf{N}}x,y\left[ x\mathbf R\sqsubset_{\mathscr P_{%
\mathbf{D}}}y\mathbf S\Leftrightarrow x\mathbf R^{\prime} \sqsubset_{\mathscr P_{\mathbf{D}%
}}yS\right],
\]
which shows that there as many pairs of binary relations $R$ and $S$ that have $\mathscr P_{\mathbf{D}}$, such that the ratio 
$\mathfrak{R}\left[ R,S\right]$ is the same as $\mathfrak{R}\left[ \mathbf R,\mathbf S\right]$, as binary relations that have 
$\mathscr P_{\mathbf{D}}$, and
there is no pair of binary relations $R^{\prime}$ and $S^{\prime}$ that have $\mathscr P_{\mathbf{D}}$
such that the ratio $\mathfrak{R}\left[ R^{\prime},S^{\prime}\right]$ is
distinct from all ratios $\mathfrak{R}\left[ \mathbf R,S\right]$, where  $S$ is a binary relation that has $\mathscr P_{\mathbf{D}}$.

Insofar as an analogous argument also applies, \textit{mutatis mutandis},
to the second of the three disjoints forming the right-hand side of EP, and the third of these disjoints only concerns ordered pairs of binary relations whose
first element is the identity relation and makes all ratios whose denominator is such relation identical, this is enough to
conclude that the
cardinality of the set $\{\mathfrak{R}\left[ R,S\right] :\mathscr M_{\mathbf{%
D}}R\wedge \mathscr P_{\mathbf{D}}S\}$ is the same as that of 
$\left\{ R:\mathscr M_{%
\mathbf{D}}R\right\}$, namely  $2^{\aleph _{0}}$, as was to be
proved. 

Hence, defining real numbers as $\mathfrak{R}$-\textit{abstracta} over a
domain of magnitudes entails the existence of at least $2^{\aleph _{0}}$
such numbers, since there is just one real number for any equivalence class
induced by the
right-hand side of EP or EP$^*$ on a single domain of magnitudes.
To prove that there are \textit{just} $2^{\aleph _{0}}$, recall that, as observed in \S ~\ref{RN} above, from EP or EP$^*$ it
follows that, if there were several such domains, for any $\mathfrak{R}$-%
\textit{abstractum} on one of those domains, there would just be one $\mathfrak{R}$-\textit{%
abstractum} on the others that is identical with it, which entails
that, if there were several domains of magnitudes, the $\mathfrak{R}$-%
\textit{abstracta} on one of them would be just the same objects as the $%
\mathfrak{R}$-\textit{abstracta} on the other.

\section{Logicality and Arithmeticity\label{logicality}}

\noindent Given all the previous considerations, we can finally
tackle two major questions concerning the definition of reals in 
FMR or FMR$'$: 
\textit{i})~Is there a strong enough sense in which this definition is
logical? \textit{ii})~Is this definition independent of natural
numbers and their theory? Insofar as it seems difficult to imagine a
consistent definition which is closest to Frege's envisaged one, the answer to these questions is
relevant to assess Frege's achievements as well: Was Frege's plan
for defining real numbers as ratios of magnitudes compatible with a
logicist program, inconsistency apart? Was it in line with the basic idea that real and
natural numbers are essentially independent objects? The question is not only a historical one, it has also a contemporary
philosophical relevance: Is a neologicist program concerning real analysis,
making it both logical and independent of the arithmetic of natural
numbers, envisageable along Frege's original indications? In what follows we
will suggest a negative answer to all these questions\footnote{Simons (1987,
\S 7) has stressed the crucial differences between
Frege's logicism for natural numbers and his views on real ones. Without undermining
his arguments---which take however for granted the usual reading
of Frege's logicism for natural numbers, which we rather take as
questionable under many respects: See Panza (2018) and Panza (FC2)---we
follow another strategy here: we frontally attack the idea that Frege's envisaged 
definition of real numbers might be taken as logical in any substantial sense.}.

All of them boil down to two issues: whether FMR or FMR$'$ can 
be taken to be logical systems, independent of a previous definition of natural numbers (likely got through FA); whether an existence proof of nonempty domains of magnitudes
and of real numbers as defined in FMR or FMR$'$ is compatible with
the logicality and arithmeticity of these definitions.

\subsection{About the Definition of Domains of Magnitudes\label{DDM}}
Insofar as FMR and FMR$'$ are obtained by adding 
some new axioms to L$_2$PCA$^{2}_{\Delta^{1}_{0}}$, we will begin by 
investigating 
whether this latter system
is genuinely logical 
and independent from the natural numbers. Both issues also apply to our definition
of domains of magnitudes within it.

Likely, no one would question
its independence from the natural
numbers. The considerations advanced in \S ~\ref{Sect.LE} seem, moreover, to support its logicality\footnote{But
see footnote (\ref{Quine}), on this matter.}. Still, admitting that L$_2$PCA$^{2}_{\Delta^{1}_{0}}$ indeed has these features
is not enough for concluding that our definition of domains of
magnitudes is, in turn, independent of natural numbers and logical in some
more significant sense than the simple and quite weak one of being formulated within a logical system. 
There are two concerns, here. 

The first is that, even in a logical system, it seems possible to
define items whose logical nature is suspect. Panza (2018) and Panza (FC2) already raised the question in relation both to natural numbers and magnitudes, as originally
defined by Frege as appropriate extensions.
Surely, according to our reformulation of Frege's definition, magnitudes are no more extensions, but rather binary first-level
relations. Still, apart from the identity relation, the relations forming such domains are not identified as particular relations somehow precisely defined;
they are rather characterized as possible places in whatever system
exemplifying a certain structure. This makes this definition define domains of magnitudes, but not magnitudes
as such, which is perfectly in line with Frege's remark quoted in \S~\ref{FS}. Hence,
all that the definition fixes is
the structure of a domain of magnitudes,
not its content, which is to be given independently of it.

As such, this
might even 
be taken as an argument for its logicality, if, \textit{contra} Frege, it is admitted that 
logic has no content. But it
makes the second concern crucial. As already claimed, whereas an existence
proof of nonempty domains of magnitudes cannot be provided within
L$_2$PCA$^{2}_{\Delta^{1}_{0}}$ (as well as within FMR and, \textit{a fortiori}, FMR$'$), it 
is indispensable for making our definition of real numbers sensible. 
So, proving, necessarily outside these systems, the existence of a nonempty domain of
magnitudes is an essential part of this
very definition (even if this is not required to formulate the definition of domains of magnitudes themselves), not only of a model-theoretical enquiry on
it. This makes both the logicality and the arithmeticity of the definition
crucially and questionably depend on the means,
external to FMR and FMR$'$, needed for conducting such an existence proof.
Two major issues arise. 

The first is that it seems plausible to require that 
a genuinely logical definition not be in need of an 
external existence proof of the items it
defines---or of other associated ones. Since such a definition should  purportedly
ensure the existence of these items 
by merely showing that, if there
were none, some logical, or innocent enough, truths would not be true after all. This
is indeed allegedly the case of the neologicist existence proof of natural
numbers\footnote{
To see it, consider the argument proving the existence of $0$
that we have detailed at the beginning of \S ~\ref{existence} above. The
concept $\left[ x:x\neq x\right] $ exists not only by predicative, but also
by logical comprehension, as well as logic is enough to get that $\left[
x:x\neq x\right] \approx \left[ x:x\neq x\right] $, and HP is so
for getting that $0=0$, which could not be true if $0$ did not exist. The
same pattern allows to prove the existence of each natural number,
provided that comprehension be extended to formulas involving the operator `$%
\#$'. So, $1$ is proved to exist, for example, since, if it did not, it
would be false that $1=1$, which follows, in agreement with HP,
from `$\left[ x:x=0\right] \approx \left[ x:x=0\right] $', which follows, in
turn, by logic, from the existence of the concept $\left[ x:x=0\right] $,
ensured by comprehension applied to the formula `$x=\#\left[ x:x\neq x\right]
$'. On the other side, the existence of the totality of numeral numbers is proved by proving, by HP and 
impredictaive comprehension, the successor axioms, which would be false, if these numbers did not exist.}---which we mimicked in the existence proof of the $\mathsf{E}_{\oslash }$-\textit{abstracta}, in \S\ \ref{Sect.KS} above.
One might argue that 
this is too demanding.
However, a distinction should be drawn between definitions that
are 
logical in this demanding sense, and others that are not or
cannot be so. This would be enough for concluding that neither our definition of real numbers in FMR or FMR$'$, nor that of domains of magnitudes
in L$_2$PCA$^{2}_{\Delta^{1}_{0}}$, can pretend to be logical in
the same sense in which neologicists claim their definition of natural
numbers is\footnote{In fact, neologicists usually take their definition of natural numbers to be
analytic, though not logical. Still, we made clear from the very beginning why 
we do not endorse this distinction
here---see \S~1 above. Let us notice, however, that the argument just advanced is emblematic
of the reason we advanced to justify our attitude.
Since, if this distinction were admitted, this argument should allow to conclude also
that our definition of real numbers is no more analytic in the same sense in which
neologicists take their definition of natural numbers to be so.
\label{ft.NDANL}}.

This leads to the second issue: once admitted that the neologicist's proof-pattern does not
apply, and that this prevents our Fregean definitions of domains of magnitudes and real numbers to be
logical in the above demanding sense, the question arises whether
these definitions might nevertheless be deemed logical in some
less demanding sense\footnote{We leave here apart  the question of whether the
neologicist definition of natural numbers or our definition
of the $\mathsf{E}_{\oslash }$-\textit{abstracta} are actually logical or analytical. In
Panza (2016), pp. 420-423, the point is made that the former definition might be deemed so in
a quite peculiar sense, quite different than those that are current in the discussion on logicism and neologism,
and because of a completely different argument than the neologicist's. There is no need to come back on 
this point, here. It is only important
to observe that it does
in no way apply to our definition of domains of magnitudes and real numbers.}.
The question seems  to have different answers according to whether  
it concerns the former definition or the latter, and whether domains of magnitudes
are regarded as such or as tools for defining real numbers.
If we look at the definition of domains of magnitudes as such, and admit that  
L$_2$PCA$^{2}_{\Delta^{1}_{0}}$ be a genuine logical system,
it is hard to find any other reason than that raised above to
deny its logicality. But if we look at
domains of magnitudes as tools for defining real numbers,
the situation changes. Insofar as proving
the existence of nonempty such domains is essential
for enabling them to play this role, the question becomes whether 
the proof can be so shaped as to make it logical, and, then, part
of a logical definition of these numbers. This is, then, the question we have to
tackle, now.

Above we explored two different strategies for conducting this proof: an
inflationary and a non-inflationary one. In what follows, we will expand on them by considering how they score with respect to the issue of logicality. Insofar as the question of the logicality of our definition of real numbers
has multiple interconnections with that of its arithmeticity, we will also consider in the meantime
whether these strategies can make this definition non-arithmetical.

\subsubsection{About the Inflationary Strategy\label{Sect.AIS}}

\noindent A proof following the inflationary
strategy may be deemed non-logical
just because of
its inflationary nature. The reason is obvious:
insofar as no proof by countable induction is possible here, such a proof cannot but grant that 
the abstraction principle introducing the continuously 
many objects it concerns (FP, in our case) \textit{eo ipso} entails the existence of these objects, which appears to
be incompatible with its being logical, and, then, part of a logical proof.

One could object that the argument merely
points out that the proof is not logical, because it requires  means other than countable induction to prove the existence of continuously many objects, which is unfair,
at best. After all, real numbers
must be continuously many, so that accepting
this argument would amount to principledly excluding the possibility of a logical definition of real numbers ensuring their existence.  The objection is not convincing. It is entirely possible that no such definition be logical. Still, if a definition of these numbers is offered with the aim of being so, it should, at least, avoid requiring an existence proof for continuously many objects  other than the reals. This would leave room for arguing that proving the existence of continuously many objects is not part of its job, but should be left for further meta-theoretical considerations. The point is, then, that the inflationary strategy is not suitable for entering a logical definition of real numbers because it requires an existence proof of continuously many objects other than the reals: a proof which cannot but appeal to independent resources from those involved in the definition itself. 

Another essential feature of the inflationary strategy is also relevant for the present discussion: 
its delivering
an arithmetical copy of the additive group of the
reals as a condition for making their definition appropriate.
This makes clear both its arithmetical nature, and its essential
mathematical circularity. Let us consider the two allegations in turn. 

To reply to the arithmeticity allegation, it is not enough to argue that proving the existence of a nonempty domain of magnitudes arithmetically does not make a definition of real numbers as ratios of magnitudes arithmetical. The fact that an existence proof of such a domain is an indispensable part of the definition immediately entails, indeed, that this definition can be deemed non-arithmetical only in the presence of an existence proof of a non-arithmetical nonempty domain of magnitudes. Since, if one could only prove the existence of arithmetical such domains, defining real numbers as ratios on them would make them arithmetical items, after all\footnote{The point might be softened by observing that our Fregean definition of domains of magnitudes differs from other possibly arithmetical ones for not appealing  to any specific property of the objects on which the relevant binary relations are defined. 
To better see this, we can compare this definition to one mimicking Dedekind's definition by cuts in terms of binary relations (we thank Andrew Moshier for his suggestion). Let $\left( \mathbf{O},<_{\mathbf{O}}\right)$ be a totally ordered set without endpoints, whose elements count as objects. By adopting a third-order logical system with
third-order predicative comprehension, the following explicit definition can be provided, where the index `$\mathbf{O}$' restricts the quantifiers to binary $\mathbf{O}$-relations (i.e. binary relations among the elements of $\mathbf{O}$) and to these very elements, respectively:
\begin{equation*}
\forall_{\mathbf{O}} R\left\{ \mathscr C_{\mathbf{O}}R\Leftrightarrow 
\left[\begin{array}{l}
\forall _{\mathbf{O}}x,y\left[ xRy\Rightarrow x<_{\mathbf{O}}y\right] \wedge
\medskip  \\ 
\forall _{\mathbf{O}}x,y,z,w\left[ \left( xRy\wedge z<_{\mathbf{O}}x\wedge
y<_{\mathbf{O}}w\right) \Rightarrow zRw\right] \wedge \medskip  \\ 
\forall _{\mathbf{O}}x,y\left[ xRy\Rightarrow \exists _{\mathbf{O}}z,w\left[
zRw\wedge x<_{\mathbf{O}}z\wedge w<_{\mathbf{O}}y\right] \right] \wedge
\medskip  \\ 
\forall _{\mathbf{O}}x,y\left[ x<_{\mathbf{O}}y\Rightarrow \left(\exists _{\mathbf{%
O}}z\left[ xRz\right] \vee \exists _{\mathbf{O}}w\left[ wRy\right]\right) \right]
\wedge \medskip  \\ 
\forall _{\mathbf{O}}x,y,z,w\left[ \left( xRy\wedge zRw\right) \Rightarrow
xRw\right] 
\end{array}%
\right]
\right\} .
\end{equation*}
The third conjunct of the right-hand side entails that no binary $\mathbf{O}$-relation has the property $\mathscr C_{\mathbf{O}}$ if $%
\left( \mathbf{O},<_{\mathbf{O}}\right) $ is not dense.\ Let us suppose that it be so. Call a relation `$\mathscr C_{\mathbf{O}}$-relation' if it
is a binary relation having $\mathscr C_{\mathbf{O}}$. We can say that any $\mathscr C_{\mathbf{O}}$-relation defines a cut on $\left( \mathbf{O},<_{%
\mathbf{O}}\right) $. The collection of the $\mathscr C_{\mathbf{O}}$-relations does not form a domain of
magnitudes, in the sense established above, since the $\mathscr C_{\mathbf{O}}$-relations are not permutations. Still, we might weaken Frege's requirement on domains of magnitudes, and take
such domains as constituted by totally ordered, dense and Dedekind-complete groups of first-level binary relations, independently of their being permutations. The $\mathscr C_{\mathbf{O}}$-relations might, then, form a domain of magnitudes if a commutative (and
associative) addition admitting a neutral element be defined on them. To this purpose, let us
suppose that an addition $+_{\mathbf{O}}$ be defined on $\left( \mathbf{O},<_{\mathbf{O}}\right)$, so as to make $\left( \mathbf{O},<_{\mathbf{O}%
},+_{\mathbf{O}}\right) $ a totally ordered, Abelian and dense additive group. We can easily define an  addition $+_{\mathscr C_{\mathbf{O}}}$ on the 
$\mathscr C_{\mathbf{O}}$-relations, by stating that
\begin{equation*}
\forall_{\mathscr C_{\mathbf{O}}}R,S\forall _{\mathbf{O}}x,y\left[ x\left( R+_{%
\mathscr C_{\mathbf{O}}}S\right) y\Leftrightarrow \forall _{\mathbf{O%
}}z,w,
v,u \left[\left( zRv\wedge wSu\right) %
\Rightarrow 
\left(x<_{\mathbf{O}}z+w \wedge v+u <_{\mathbf{O}}y\right)\right]\right],
\end{equation*}
where the index `$\mathscr C_{\mathbf{O}}$' to the first universal quantifier restricts it to these relations. The $\mathscr C_{\mathbf{O}}$-relation $\mathbf{Z}_{\mathbf{O}}$ defined by
\begin{equation*}
\forall _{\mathbf{O}}x,y\left[ x\mathbf{Z}_{\mathbf{O}}y\Leftrightarrow x<_{%
\mathbf{O}}0_{\mathbf{O}}<_{\mathbf{O}}y\right] 
\end{equation*}%
is the neutral element of $+_{\mathscr C_{\mathbf{O}}}$, and another $\mathscr C_{\mathbf{O}}$-relation $\mathbf{R}$, is deemed positive if and only if 
$$\exists_{\mathbf{O}}x,y\left[xRy \wedge 0_{\mathbf{O}} < _{\mathbf{O}} x <_{\mathbf{O}}y\right],$$ and negative otherwise. One could, then, define an order relation $\sqsubset_{\mathscr C_{\mathbf{O}}}$ on the $\mathscr C_{\mathbf{O}}$-relations by stating that 
\begin{equation*}
\forall_{\mathscr C_{\mathbf{O}}}R,S\left[ R\sqsubset_{\mathscr C_{\mathbf{O}}} S\Leftrightarrow \exists _{\mathscr C_{\mathbf{O}}^{+}}T\left[ R+_{%
\mathscr C_{\mathbf{O}}}T=S\right] \right],
\end{equation*}%
where the index `$\mathscr C_{\mathbf{O}}^{+}$' restricts the existential quantifier 
to positive $\mathscr C_{\mathbf{O}}$-relations. These would form a totally ordered, dense and Dedekind-complete group under $+_{\mathscr %
C_{\mathbf{O}}}$ and $\sqsubset_{\mathscr C_{\mathbf{O}}}$, and, then, a domain of magnitudes, in the previous weakened sense. One might, then, define real numbers as ratios on such a group. Still, so defined, real numbers would be, structurally speaking, nothing more than ratios on cuts-relations on the additive group of the rational numbers, and this would make them intrinsically arithmetical items, in a much stronger sense than the real numbers
defined in FMR or FMR$'$, under the condition that the only nonempty domains of magnitudes whose existence can be proved were arithmetical.}.

To reply to the circularity allegation, one should argue that the
copy of the additive group of the reals is
just a copy, since, though structurally coincident with real numbers, its 
elements intrinsically differ from them.\ We can imagine
Frege advancing this argument.\ But we can hardly follow him in  this
without making any working mathematician sarcastically smile.

Let us recap. The inflationary strategy suffers from two problems: in
absence of a further existence proof of non-arithmetical 
domains of magnitudes, it makes real numbers themselves arithmetical objects, after all; 
it requires a preliminary structural definition of the real numbers, in order
to make the planned definition of these same numbers suitable\footnote{Both problems
also arise if the inflationary strategy is implemented by 
$\mathsf{E}_{\oslash }$-\textit{abstracta}. As for the
first, this is obvious. As for the second, notice that these
\textit{abstracta} could enter the existence proof only
because of their features that make them structurally
coincide with the natural numbers.\label{ft.SIEASNN}}.

\subsubsection{About the Non-Inflationary Strategy}
At least four reasons might be advanced to argue that, in the light of the
existence proof in \S\ \ref{Sect.KS} above, our Fregean definition
of real numbers is hardly both logical and
non-arithmetical: \textit{i})~that proof is based on 
the symmetric group $\mathbf{\Sigma }_{\mathbb{N}}$ on the
natural numbers; \textit{ii})~it allows to conclude that a nonempty domain
of magnitudes exists only if it is admitted that the symmetric group on an
infinite countable set exists if this set exists; \textit{iii})~it
essentially appeals to the additive group $\left( \mathbb{R},+\right) $ of
the real numbers themselves; \textit{iv})~it appeals to Zorn's lemma, and is, then, not constructive.

The first reason cannot be dismissed by merely observing 
that, in our reconstruction of the proof, 
$\mathbf{\Sigma }_{\mathbb{N}}$ has been replaced by the
symmetric group on the set of $\mathsf{E}_{\oslash }$-\textit{abstracta}.
Since, once $\mathbf{\Sigma }_{\mathbb{N}}$ is replaced in 
Karras and Solitar's proof by any other
symmetric group over any infinite countable set, the problem becomes that of
justifying the existence of this set, and even its cardinality, by no appeal to
$\mathbb{N}$ itself. That $\mathbf{\Sigma }_{%
\mathbb{N}}$ is isomorphic with the symmetric group over any infinite
countable set is, indeed, simply because any such set can be put into a bijection
with $\mathbb{N}$, so that its elements can be taken either to
count as natural numbers or, at least, to be encoded by them.

The second reason is similar to one discussed above as for the inflationary strategy: on what
logical ground can we argue that the existence of a countable
infinite set entails the existence of the symmetric group over this
set---or, more in general, of an uncountable set somehow generated by it by
considering at once some totality of properties, relations or functions
defined on the elements of this set? The fact that Frege himself suggests making a
similar admission, in order to prove the existence of a nonempty domain of
magnitudes, in no way makes it
logically licensed. Rather, it
seems to show that the very proof Frege suggested would have actually been not logical.

The third and fourth reasons are by far more delicate, and somehow interconnected.
Since, if a constructive proof of Karrass and Solitar's
theorem were available, one could hope to rely on it in order to constructively define 
a totally ordered, dense and Dedekind-complete subgroup of $\mathbf{\Sigma }_{\mathbb{N}}$,
without recurring to $\left( \mathbb{R},+\right) $.

To dismiss the third reason, and the circularity allegation that goes with it,
one might replace Karrass and Solitar's theorem with a more
general result not involving $\left( \mathbb{R},+\right) $.
A natural candidate is a result by de Bruijn (1964, p.~594), 
according to which, for any infinite cardinal $\kappa $, every Abelian
group of order $2^{\kappa }$ can be \textquotedblleft embedded
into\textquotedblright\ the symmetric group of a set of cardinality $\kappa $\footnote{This theorem was firstly published 
one year later than
Karrass and Solitar's (de Bruijn 1957, pp.
560-61 and 566), but it was then erroneously proved. The error lied with
a lemma proved by erroneously
supposing that a certain arbitrary group could be non-Abelian. The proof was later
corrected and made independent of this lemma---and in fact simplified.}. Still, the basic idea of de Bruijn's proof
is not so different from Karrass and Solitar's and makes this proof also depend on the axiom of Choice,
though avoiding appealing  to vector spaces. The fact that the theorem does not specifically involve $\left( \mathbb{R},+\right)$ is,
moreover, far from being an advantage in our perspective. Since  it makes this theorem unable 
to provide a ground for the required existence proof.
For the purpose of this latter  
proof is establishing  that $\mathbf{\Sigma }_{%
\mathbb{N}}$ (or, more generally, the group of permutations on a countable
set) actually includes a subgroup complying with the
relevant structure, while this theorem merely ensures that,
if there is such a group, then it can be embedded into $\mathbf{\Sigma }_{%
\mathbb{N}}$, and can be regarded
as a subgroup of it. This makes, of course, de Bruijn's theorem immediately
entail that $\left( \mathbb{R},+\right) $ can be embedded into $\mathbf{%
\Sigma }_{\mathbb{N}}$.\ This cannot but make the circularity even more evident, since 
it is only the existence of $\mathbb{R}$ that can ensure that
a subgroup of $\mathbf{\Sigma }_{\mathbb{N}}$ complying with the relevant
structure exists. In order to solve
the issue, one should prove that $\mathbf{\Sigma }_{\mathbb{N}}$
includes a totally ordered, dense and Dedekind-complete 
subgroup, without assuming the existence of this group. To the best of our knowledge, this has not yet been done.

This does not mean, of course, that this result, or any other entailing it, has
not actually been proved or, even less so, that
this cannot be done. The fourth reason suggests, however, that the relevant
question is not whether this has been or might be done, but, rather,
whether this can be done constructively, i.e. without appealing to a form of the Axiom of Choice,
which might hardly be taken as a logical principle. When put in a clear
mathematical form, the question is whether it is provable in ZF alone
(or in some other appropriate setting that neither presupposes nor entails the Axiom of Choice)
that $\mathbf{\Sigma }_{%
\mathbb{N}}$ contains a totally ordered, dense and Dedekind-complete
subgroup, and whether, moreover, this can be done without 
assuming the existence of $\left( 
\mathbb{R},+\right) $. To the best
of our knowledge, this question also has not been answered yet.

To begin enquiring about it, one might wonder whether Karrass and Solitar's
proof can be freed from Zorn's lemma or any equivalent assumption. 
Such an
assumption enters the proof to ensure that any vector space has a basis---i.e. that 
such a basis exists though it cannot be
constructively displayed. 
This makes it relevant to observe that Blass (1984)
proves that the assumption that any vector space has a basis is
(ZF-)equivalent to the Axiom of Choice. This is still not enough to ensure that the
appeal to a form of the Axiom of Choice cannot be
avoided in Karrass and Solitar's proof, and that this proof is, then, both non-constructive and intrinsically
dependent on such an axiom. Since what is required for this proof to work is not, properly,
that any vector space has a basis, but rather that this is so for the two relevant
such spaces, i.e.  $\prod\limits_{i\in \mathbf{I}}\left( \mathbb{Q}%
,+\right) $ and $\left( \mathbb{R},+\right) $. The issue becomes, then,
whether one can prove that these very vector spaces have a basis without appealing
to a form of the Axiom of Choice or to any other non-constructive means. 
To the best of our knowledge, anew, this is still unknown. 

Still, even if this could not be done, it would not follow that Karrass and
Solitar's theorem cannot be proved without appealing to a form of this axiom. 
The only occurrence of Zorn's lemma in the previous proof is,
indeed, in its very last step, which is the only one involving
vector spaces. It is, then, natural to wonder whether the theorem could be proved by avoiding
this step (and, then, presumably any
reference to vector spaces), by replacing it with another step not depending on the Axiom of
Choice or some other non-constructive assumptions.

We could imagine two scenarios. In
the first, the question is whether ZF alone is capable of proving
Karrass and Solitar's theorem: to this extent, either it is, or the theorem
is undecidable there. In the second scenario, the question is whether
ZF augmented with some axioms incompatible with the Axiom of Choice, such as
the Axiom of Determinateness,\footnote{%
First introduced by Mycielski \& Steinhaus (1962), this axiom asserts that
\textquotedblleft certain infinite, deterministic 2-person games with
complete information [\ldots ] are determinate, i.e., that one of the
players has a winning strategy\textquotedblright .\ See also Herrlich (2006, p. 151), 
which also provides a proof of the incompatibility between
the Axiom of Choice and the Axiom of Determinateness.} is capable of
proving this theorem.\ To the best of our knowledge, these issues
have also not been settled yet\footnote{%
There still might be clues on the second issue: while the
Axiom of Choice entails that $\mathbb{R}$, as a vector
space over $\left( \mathbb{Q},+,\cdot \right) $, has bases, the Axiom of
Determinateness entails that it has not (Herrlich 2006, Theorem 4.44 and
Corollary 7.20)---and one might even guess that it be the same for an
infinite product of copies of $\left( \mathbb{Q},+\right) $.} .

The conclusion to be drawn from all these remarks cannot be but prudent.
Still, it can certainly no more suggest that Karrass and Solitar's theorem provide a basis to
argue that our Fregean definition of real numbers is logical. Since, circularity issues aside, 
arguing that this theorem can enter a
non-inflationary existence proof for a nonempty domain of magnitudes,
suitable for making our definition logical, would be
quite premature, at best. And it would be even more so 
to argue that
a more general theorem, asserting that $\mathbf{\Sigma }_{\mathbb{N}}$
includes an appropriate group identified without
appealing to $\mathbb{R}$, can enter such a non-inflationary existence
proof.

\subsection{Might the Existence Proof be Avoided?}

\noindent The previous considerations suggest that there is
no way to prove the existence of
nonempty domains of magnitudes without  wiping out  both
the logicality and non-arithmeticity of our Fregean definition of real numbers.
In light of this conclusion, one might suggest
changing the rules of the game. Even if there is no way
to prove, by (higher-order) logic, suitable existentially innocent
abstraction principles and appropriate algebraic and/or set-theoretical
constructive arguments, that nonempty domains of magnitudes exist,
still we know they do. For we can show or prove it, for example, by empirically-tied 
geometric or mechanical considerations; or, as above, by
trusting non-constructive set-theoretical principles; or by assuming that
natural numbers exist and granting the previous abstraction principles an existential import.
Hence, one could argue, defining the real numbers 
\textit{\'a la} Frege within FMR or FMR$'$, even with no
existential proof, allows one to show that ratios on any externally
given domain of magnitudes, whether intrinsically arithmetical or
not, are real numbers. Since, as we have already seen, taking a real number to be a ratio on distinct 
domains of magnitudes is nothing but describing the same object in
different ways---or giving different names to it.

The problem with this move is that applying our definition to whatever
externally given domain of magnitudes would certainly warrant that the
ratios on it are real numbers, but not that real numbers are
intrinsically such ratios, let alone that they are
non-arithmetical items. If we reasoned this way, we would
do nothing essentially different from appealing to a representation theorem to 
draw the conclusion that real numbers
measure the magnitudes in the relevant domains, in the spirit of the measurement theory\footnote{%
About the tension between considering applications of real numbers in
agreement with the measurement theory and taking them to be ratios on
domains of magnitudes, see the Hale-Batitsky discussion in Hale (2000, 2002) and
Batitsky (2002). On this matter see also Panza \& Sereni (2019, pp.\ 122-123
and 126-130).}. In both cases, all we do is recognize that some
externally given systems (arithmetical or not) comply with some fixed
structural conditions. The fact that these structural conditions are fixed by
our definition in FMR or FMR$'$, or by recurring to algebraic
axioms as those of a totally-ordered complete Abelian field (as usually
supposed in measurement theory), or, again, by alternative
definitions (as Cantor's and Dedekind's, by Cauchy's sequence and
cuts on rationals, respectively, or even as the one
grounded on FP)
makes no essential difference on this matter.

In the eyes of a
Frege partisan, there would be a crucial difference only if the existence proof were deemed an
essential, though supplementary, part of the definition itself, as we did
above. Since this would make the numbers so defined intrinsically
ratios on domains of magnitudes, and
their application to measurement
\textquotedblleft built into\textquotedblright\ their nature and/or their
very definition, as required by the application constraint (Wright, 2000, p.~325). 
In this respect, the previous remarks
on the arithmeticity and logicality of our definition in FMR or FMR$'$
should be intended to suggest that compliance with this
constraint is incompatible not only with offering a logical definition of
real numbers, as already argued in Panza \& Sereni (2019), but
also with defining them non-arithmetically,  
despite Frege's
adhesion to the same constraint as the main source of his quest for a non-arithmetical definition
of these numbers.

\subsection{About Euclid's Principle}
Up to now, we have only considered the existence proof of nonempty
domains of magnitudes.\ Still, the indispensability and the nature
of this proof are not the only reasons suggesting that our Fregean
definition of real numbers is neither logical nor non-arithmetical. Since,
once domains of magnitudes have been defined and somehow proved to exist,
the question remains open of defining real numbers as ratios on them. 
In our setting, this is done by means of EP$^*$\footnote{Possibly with the help of
(CA$^{2}_{\Sigma^{2}_{1}}$), if an explicit definition  like (\ref{5.1}) is required. 
We do not want to enter here a discussion about the logicality of (the different sorts of)
comprehension. We merely remark that the high impredicativity of 
(CA$^{2}_{\Sigma^{2}_{1}}$) might make many doubt not only its
logicality, but also its licitness. Who doubts both has no other choice but rejecting
definition (\ref{5.1}) and rest content with (\ref{5.1}$'$).
Who doubts only the former can either admit definition (\ref{5.1}), but take it as non-logical,
or rejecting, again, this definition in favor of (\ref{5.1}$'$).}.
The questions are, then, obvious: is this principle logical? Is it actually independent 
of natural numbers?

Let us start from the latter. Undoubtedly, EP depends on natural numbers, as it explicitly quantifies over them.
Still, this principle
also depends on the use of other linguistic means foreign to 
L$_2$PCA$^{2}_{\Delta^{1}_{0}}$, by
involving a piece of informal language allowing for
predicate variables like `$xR$'. 
The two features are connected, since the
quantification over natural numbers just operates on the individual variables
occurring in 
these predicate variables. Replacing 
EP with EP$^*$ allows avoiding
both the quantification over natural
numbers and the appeal to informal language at once.
Surely, EP$^*$ involves
no second-order predicate constant supposedly designating the
property of being a natural number. Still, this does not
ensure, yet, its independence from natural numbers, since it is far from clear that  the trick used to
avoid the quantification on these numbers is actually independent of them. What might make one 
suspect it is not that 
the right-hand side of (\ref{eq.A}) is just an appropriate third-order
rephrasing of the right-hand side of HP. Hence, if it were admitted that, 
when applied to finite concepts, HP
is intrinsically
inherent to natural numbers---not only because the objects not
complying with it are not natural numbers, but also because its
assumption \textit{ipso facto} brings these latter about---, one should
infer that, in spite of appearance, also
EP$^*$ depends on these numbers. As a matter of fact, this is a strong
assumption, but one that can be made in a Fregean vein, and which might bring, then, \textit{ipso facto}---that is, independently of any consideration on the existence proof of nonempty domain of magnitudes---,
to the conclusion that  
our
definition of real numbers, whether in FMR or FMR$'$, is essentially 
arithmetical.

Someone admitting this assumption might still argue against this
conclusion by observing that (\ref{eq.A})
essentially differs
from HP for being a (metalinguistic) abbreviation stipulation, rather than an
axiom providing an implicit definition of a functional constant. This is
enough, one might continue, to make EP$^*$
appeal to no variable ranging on objects that might count as the
natural numbers. This is unquestionably so. However, any instance of 
`$_{\left( \mathscr X,\mathscr X^{\prime }\right) }\mathcal{E}\left(
 R, T, R^{\prime }, T^{\prime }\right)$' asserts that a certain 
first-level binary relation
is the same multiple of another such relation over a certain positive class
as a third such relation of a fourth one, over
the same or another positive class. This is, in turn, the same as
asserting that the iterations of the composition operation on such a
relation within the former class are into a bijection with the iterations of
composition on such a relation within the latter class. If
this is not the same as making natural numbers enter into play, it
is, at least, the same as making
the equinumerosity
relation so. Hence, if EP$^*$
is not dependent on natural numbers, it seems to be, at least, dependent on
counting. There is no easy way to settle whether this is enough to make EP$^*$ an arithmetical principle. Here,
we just observe that this makes our Fregean definition of
real numbers, whether in FMR or FMR$'$, dependent 
on an essential ingredient of any Fregean
definition of natural numbers. Even if this, as such, does not make our definition
arithmetical, it is plausibly enough for making it much more related to natural
numbers than Frege might have desired his definition be.

Let us come, now, to the first question: can EP$^*$ be deemed logical?
A simple way to tackle the question might be that of choosing between two quite natural options:
either any abstraction principle is logical if it is stated through a logical
language, or it is so only because
of the peculiar nature of its
right-hand side. In the former case, EP$^*$ is logical if  L$_2$PCA$^{2}_{\Delta^{1}_{0}}$ is so.
In the latter case, EP$^*$ 
cannot be logical on the same grounds on which HP or a consistent
version of BLV might be so. If this simple alternative
is rejected, if only for argument's sake,
what criterion might be provided to distinguish logical
abstraction principles stated in a logical language 
from non-logical ones?
Consistency is surely not enough.
But, then, what? We cannot dwell on this
issue here. We simply contend that the burden of the proof seems to be on anyone arguing that EP$^*$ is logical,
despite its being essentially akin to Euclid's definition of proportionality,
which has been considered for centuries as the cornerstone of the most
fundamental mathematical theory on which classical geometry was crucially
grounded.

\section{Concluding Remarks}

\noindent Though some of them are certainly far from knock down ones, we think we advanced enough arguments in favor of the claim that our rendering of Frege's envisaged definition of real numbers is
neither logical nor non-arithmetical. As our rendering is arguably the
closest possible to it, this conclusion questions the possibility that Frege's own definition could be achieved logically and non-arithmetically. It remains to establish whether this was actually
Frege's intent.

That Frege was aiming at a logical definition of real numbers as his main goal
for the foundation of real analysis might be questioned for several
reasons\footnote{%
Some are advanced in Panza (FC2).}. One of them might be the
following. 

From our reconstruction, it seems to emerge that arguing for the logicality of a definition of real numbers following Frege's
indications requires arguing that FP,
or any akin
principle, is both logical and capable of delivering continuously many
objects without the assistance of
any independent existence proof. But if so, then  
FP would also be enough for delivering real
numbers, if not as logical objects, at least as objects defined in a
logical setting. But, then, why did Frege
venture
himself in a so entangled definition whose logical nature is as
suspect as that of real numbers as ratios of magnitudes? 

Possibly, far from considering logicality as his ultimate aim, he overall wanted
to link real analysis to a general theory of magnitudes. This has been argued for in
Panza \& Sereni (2019). Or, possibly, he merely wanted to distinguish real from natural
numbers, making the former essentially independent of the latter, for their
being objects of an essentially different kind. Though the two possibilities are not incompatible with each other, our 
conclusion might be taken as a piece of evidence that he could
not have reached this second aim by following the route envisaged in the
\textit{Grundgesetze}. The first aim remains, which is
certainly paramount from a purely mathematical perspective. If we admit that
this was, after all, his prominent goal, then our rendering of his definition
might be taken as an indication of a simple way to accomplish it.\\

\noindent \textbf{Acknowledgments:} We would like to thank Mirna D\v{z}amonja, Alain Genestier, 
Paolo Mancosu, Andrew Moshier, Jamie Tappenden, an anonymous reviewer for \textit{The Review of Symbolic Logic}, and the 
audiences of the conferences previous drafts of this paper were presented at, for useful comments and criticisms. We are also indebted to Mya MacRae for the linguistic revision of the paper. The work of Marco Panza has been supported by the ANR-DFG project FFIUM.

\end{document}